\documentclass[10pt]{article}%
\usepackage{amsmath}
\usepackage{amsfonts}
\usepackage{amssymb}
\usepackage{graphicx}%
\newtheorem{theorem}{Theorem}

\newtheorem{corollary}[theorem]{Corollary}

\newtheorem{definition}[theorem]{Definition}

\newtheorem{lemma}[theorem]{Lemma}

\newtheorem{proposition}[theorem]{Proposition}

\newenvironment{proof}[1][Proof]{\noindent\textbf{#1.} }{\ \rule{0.5em}{0.5em}}
\begin{document}

\title{\textbf{Being serious about non-commitment: subgame perfect equilibrium in
continuous time}}
\author{\textbf{Ivar Ekeland\thanks{University of British Columbia, mathematics and
economic departments (ekeland@math.ubc.ca)} and Ali Lazrak\thanks{Sauder
School of Business at the University of British Columbia
(ali.lazrak@sauder.ubc.ca)}.}}
\date{Preliminary version: Mars 2006}
\maketitle

\begin{abstract}
This paper characterizes differentiable subgame perfect equilibria in a
continuous time intertemporal decision optimization problem with non-constant
discounting. The equilibrium equation takes two different forms, one of which
is reminescent of the classical Hamilton-Jacobi-Bellman equation of optimal
control, but with a non-local term. We give a local existence result, and
several examples in the consumption saving problem. The analysis is then
applied to suggest that non constant discount rates generate an indeterminacy
of the steady state in the Ramsey growth model. Despite its indeterminacy, the
steady state level is robust to small deviations from constant discount rates.

\end{abstract}

\clearpage

\section{Introduction}

This paper adresses the problem of time inconsistency under non-constant
discounting. Whereas our method and result are quite general in character, we
have chosen to illustrate them in the framework of the Ramsey model of
economic growth, (1928; see for instance \cite{Barro salai}), which has also
been used as a test case by Barro \cite{Barro} and Karp \cite{Karp} in their
investigations of the subject.

In its typical formulation, the Ramsey model represents the decision maker as
maximizing:
\begin{align}
&  \max\int_{0}^{T}h(t)u\left(  c\left(  t\right)  \right)  dt+h\left(
T\right)  g\left(  k\left(  T\right)  \right)  \ \label{3}\\
&  \frac{dk}{dt}=f\left(  t,k\left(  t\right)  \right)  -c\left(  t\right)
,\ \ k\left(  0\right)  =k_{0}\label{5}%
\end{align}
where $[0,T]$ is the life span of the decision maker and where the function
$f$ maps $[0,T]\times R^{d}$ onto $R^{d}$. Here, the decision maker can be
interpreted either as an individual or as a governement. In either case, $c$
denotes the consumption of the representative individual, $u\left(  c\left(
t\right)  \right)  $ is the utility of current consumption, $k\left(
t\right)  $ is current capital and $g(k(T))$ is the utility of terminal
capital. If the decision maker is an individual, $f(t,k(t))$ represents
capital rental interest and wages, and if it is a government, it represents
production and capital depreciation.

$h:\left[  0,\infty\right]  \rightarrow R$ is the discount function. Here and
in the sequel, it will be assumed that it is continuously differentiable, with
$h\left(  0\right)  =1$, $h\left(  t\right)  \geq0$, $h^{\prime}\left(
t\right)  \leq0$ and $h\left(  t\right)  \rightarrow0$ when $t\rightarrow
\infty$. The classical case, the one considered by Ramsey and the subsequent
litterature until the pionneering work of Strotz \cite{Strotz}, is the one
when the discount rate is constant: $h\left(  t\right)  =\exp\left(  -\rho
t\right)  .$

The decision maker, be it an individual or a government, faces this
maximization problem at time $0$, and decides on an overall solution,
$t\rightarrow\left(  \bar{c}\left(  t\right)  ,\bar{x}\left(  t\right)
\right)  $, valid for $0\leq t\leq T$. At any intermediate time $t,$ the
decision-maker, either herself at a later time if she is an individual, or
whoever is in office if it is a governement, will face a similar problem,
namely:%
\begin{align*}
&  \max\int_{t}^{T}h\left(  t-s\right)  u\left(  c\left(  s\right)  \right)
ds+h\left(  T-t\right)  g\left(  k\left(  T\right)  \right)  \ \\
&  \frac{dk}{ds}=f\left(  s,k\left(  s\right)  \right)  -c\left(  s\right)
,~~~s\geq t,\ \ k\left(  t\right)  =k_{t}%
\end{align*}
where $k_{t}$ is the existing capital at time $t$. The solution to this
problem will be some $s\rightarrow\left(  \tilde{c}\left(  s\right)
,\tilde{x}\left(  s\right)  \right)  $, valid for $t\leq s\leq T$. If this is
different from $s\rightarrow\left(  \bar{c}\left(  s\right)  ,\bar{x}\left(
s\right)  \right)  $, then the decision-maker at time $t$ is being asked to
implement a policy which, from her point of view is suboptimal. This she will
not do, unless the decision-maker at time $0$ has found a way to commit her.
If this is not the case, then the optimal policy $t\rightarrow\left(  \bar
{c}\left(  t\right)  ,\bar{x}\left(  t\right)  \right)  $ for problem
(\ref{3}) (\ref{5}) cannot be implemented. This is the problem of
\emph{time-inconsistency}, which has been studied by many authors: see
\cite{Harris-Laibson1} for a survey.

It has been known for a long time that in the case where the discount rate is
constant, so that $h\left(  t\right)  =\exp\left(  -\rho t\right)  $, time
consistency obtains: $\left(  \tilde{c},\tilde{x}\right)  =\left(  \bar
{c},\bar{x}\right)  $, so that the decision-maker at time $0$ can count on the
decision-makers at all intermediate times to implement the decisions she has
planned. The fundamental reason for which time consistency obtains is that
preference reversals due to the mere passage of time are precluded: with a
constant discount rate, relative preference between two prospective
consumption plans is unaffected by their distance into the future\footnote{To
see this, assume $0<r<s<t<T$ and consider two consumption plans $c(v)$ and
${\bar{c}}(v)$ valid for $v\geq t$. The incremental utilities ($IU$)for self
\textquotedblleft$r$" and self \textquotedblleft$s$" are related by
\begin{align*}
IU_{r}:  &  =\int_{t}^{T}e^{-\rho(v-r)}\left[  u\left(  c(v)\right)  -u\left(
{\bar{c}}(v)\right)  \right]  dv\\
&  =e^{-\rho(s-r)}\int_{t}^{T}e^{-\rho(v-s)}\left[  u\left(  c(v)\right)
-u\left(  {\bar{c}}(v)\right)  \right]  dv=:e^{-\rho(s-r)}IU_{s}%
\end{align*}
and therefore the ordinal ranking of $c$ and ${\bar{c}}$ does not change with
the mere passage of time.}. In other words, when the discount function is
exponential, relative preferences induced from the discounted utility model do
not change with time.

But why should the discount function be precisely exponential ? Experimental
evidence from psychology challenges the main consequences to be derived from
constant discount rates: see Ainslie \cite{Anslie2} and Frederick \textit{et
al} \cite{Frederick} for an overview. Relative preferences do seem to change
with time. In particular, there is robust evidence of an inclination for
imminent gratification even if accompanied by harmful delayed consequences.
This suggest a discount rate which is declining over time (see Ainslie
\cite{Anslie1} and Lowenstein and Prelec \cite{Loewenstein-Prelec}). In other
words, the discount function $h$ should be \emph{hyperbolic}, that is,
$h^{\prime}/h$ should be decreasing.

In this paper, we will deal with general discount functions:\ they need not be
hyperbolic, but they are certainly not exponential. Then time-inconsistency
obtains. We shall also assume that the decision-maker at time $0$ cannot
commit the decision-makers at later times $t>0$. This means that the solution
of problem (\ref{3}), (\ref{5}) cannot be implemented. In other words, there
is no way for the decision-maker at time $0$ to achieve what is, from her
point of view, the optimal solution of the problem, and she must turn to a
second-best policy. Defining and studying such a policy is the first aim of
this paper. The path to follow is clear. The best the decision-maker at time
$t$ can do is to guess what her successors are planning to do, and to plan her
own consumption $c\left(  t\right)  $ accordingly. In other words, we will be
looking for a subgame-perfect equilibrium of a certain game.

A second idea now comes into play: we will assume that none of the
decision-makers is sufficiently powerful to influence the global outcome. This
is very similar to perfect competition, where no agent is sufficiently
important to affect prices, and it will be formalized in the same way. In his
seminal paper \cite{Aumann}, Aumann captures that idea by considering an
exchange economy where the set of traders is the interval $\left[
0,\ 1\right]  $. An allocation then is a map $x:\left[  0,\ 1\right]
\rightarrow R_{+}^{n}$ and the total consumption of a coalition $A\subset
\left[  0,\ 1\right]  $ is the integral:%
\[
\int_{A}x\left(  t\right)  dt
\]
so that individuals, and more generally coalitions with vanishing Lebesgue
measure, have zero consumption, and therefore cannot influence prices.
However, a small coalition $\left[  t,\ t+\varepsilon\right]  $ will be able
to do so, and its weight will be roughly proportional to $\varepsilon$.

Similarly, we will consider that the set of decision-makers is the interval
$\left[  0,\ T\right]  $. At time $t$, there is a decision-maker who decides
what current consumption $c\left(  t\right)  $ shall be. As is readily seen
from the equation (\ref{5}), changing the value of $c$ at just one point in
time will not affect the trajectory. However, the decision-maker at time $t$
is allowed to form a coalition with her immediate successors, that is with all
$s\in\left[  t,\ t+\varepsilon\right]  $, and we will derive the definition of
an equilibrium strategy by letting $\varepsilon\rightarrow0$. In fact, we are
assuming that the decision-maker $t$ can commit her immediate successors (but
not, as we said before, her more distant ones), but that the commitment span
is vanishingly small.

In section \ref{sec1}, we use that idea to derive a suitable concept of
equilibrium strategy. Given a strategy $c=\sigma\left(  t,k\right)  $, a
coalition $\left[  t,\ t+\varepsilon\right]  $ will be able to perturb the
discounted utility at time $t$ by deviating unilaterally, that is, by choosing
some $c\left(  t\right)  $ different from $\sigma\left(  t,k\left(  t\right)
\right)  $; the perturbation will of course be of the first order in
$\varepsilon$. If there is no incentive for this coalition to deviate, in the
sense that this perturbation is always non-positive, and zero if and only if
$c\left(  t\right)  =\sigma\left(  t,k\left(  t\right)  \right)  $, then
$\sigma$ is an $\varepsilon$-equilibrium, in fact a subgame perfect
equilibrium. Letting $\varepsilon\rightarrow0$, we derive an appropriate
notion of equilibrium strategy in the case when individual decision makers do
not have market power.

In section \ref{sec2}, we characterize the newly defined equilibrium
strategies in terms of a value function $V\left(  t,k\right)  $. This function
is seen to satisfy two equivalent equations, (IE) and (DE), the latter being
very similar to the usual Hamilton-Jacobi-Bellman (HJB) equation of optimal
control, and reducing to (HJB) in the case when $h\left(  t\right)  $ is an
exponential. However, (DE) is not a partial differential equation: it contains
a non-local term, which makes it much more difficult to study than a
straightforward partial differential equation. We have only a local existence
result, which is stated without proof. However, in section \ref{sec3}, we
provide explicit examples in the case when the horizon is infinite, $T=\infty
$, and $f\left(  t,k\right)  $ takes the special form $r\left(  t\right)
k+w\left(  t\right)  $ (capital revenue plus wage). We also investigate the
naive strategy, where each decision-maker simply forgets that he cannot commit
his successors, and plays as if she could; we show that it is not an
equilibrium strategy, unless $u\left(  c\right)  =\ln c$.

In section \ref{sec4}, we focus on the infinite-horizon problem, with $n=1$,
and we investigate whether there is some $\bar{k}$ such that all paths
$k\left(  t\right)  $ converge to $\bar{k}$ in equilibrium. This is the
question of balanced growth, which has been much studied in the case when
$h\left(  t\right)  =\exp\left(  -\rho t\right)  $, and optimal control theory
applies; it is well known that in that case, we must have $f^{\prime}\left(
\bar{k}\right)  =r$, which effectively pins down the value of $\bar{k}.$ In
the case of general discount function, we find that $f^{\prime}\left(  \bar
{k}\right)  $ must belong to some interval, and that, \textit{ceteris
paribus}, this interval converges to the point $r$ if $h\left(  t\right)  $
converges to $\exp\left(  -\rho t\right)  $. We conclude in section \ref{sec5}.

The results obtained in sections \ref{sec4} and \ref{sec5} are very similar to
those obtained earlier by Barro \cite{Barro} and by Karp \cite{Karp}. The main
contribution of the present paper lies elsewhere, in the precise definition of
equilibrium strategies, and in their characterization through a value function
$V\left(  t,k\right)  $ which has to satisfy certain equations, reminescent of
the (HJB) equation. This allows us to carry the calculations somewhat further
than Barro or Karp, and it also opens the door to a systematic study of the
problem. The local existence result which we give is an example of what can be
obtained through our approach, and not otherwise.

The case when time is discrete, $t_{1}=0,t_{2},...,t_{n-1},t_{n}=T$, has been
investigated by many authors, for instance Strotz \cite{Strotz} , Pollak
\cite{Pollak}, Peleg and Yaari \cite{Peleg-Yaari}, Phelps and Pollak
\cite{Phelps-Pollak}, and more recently, Laibson \cite{Laibson1}. The last
decision-maker operates at time $t_{n-1}$; after he has acted, the party is
over. He is facing a plain vanilla optimization problem, and solves it. His
predecessor operates at time $t_{n-2}$. She is faced with a leader-follower
game, which she solves by integrating the strategy of her successor into her
own decision. In principle, by proceeding recursively in this way, one can go
all the way back to $t_{1}=0$, the very first decision to be made (which,
again, would not be the optimal one from the time $0$ perspective, if this
particular decision-maker could commit all her successors). If this method is
successful, it yields a subgame perfect equilibrium, and the corresponding
policy will follow through despite the lack of commitment devices. It is also
important to observe that the equilibrium policy, as in the prisoner's dilemma
game, is suboptimal relative to the outcome that can occur with a
pre-commitment technology. Using this approach (and extending it), a recent
literature has flourished showing that apparent irrationality of individuals,
even in financial markets, can be ascribed to the fact that the psychological
discount factor is not exponential; see Laibson \cite{Laibson2}, O'Donoghue
and Rabin \cite{O'Donoghue-Rabin2}, Harris and Laibson \cite{Harris-Laibson},
Krusell and Smith \cite{Krusell-Smith}, Diamond and Koszegi
\cite{Diamond-Koszegi}, Luttmer and Mariotti \cite{Luttmer-Mariotti} and others.

Unfortunately, such games typically fail to have a subgame-perfect
equilibrium. The reason is that, even if $u$ is concave with respect to $c$,
the payoff to the decision-maker at time $t_{i}$ is not concave with respect
to his own consumption $c_{i}$, because $c_{i}$ determines the capital
$k_{i+1}$ at time $t_{i+1}$, and constrains the choice of the next
decision-maker in a complicated, and certainly non-linear, way. Proceeding
recursively from $t_{n-1}$, the strategy $c_{i}=\sigma_{i}\left(
k_{i}\right)  $ at time $t_{i}$ will end up being discontinuous with respect
to $k_{i}$, which effectively kills the hope of finding a subgame-perfect
equilibrium. It is a fundamental difficulty of the discrete time model, and
various ways have been devised to get around this problem, such as adding a
public correlation device, as in Harris, Reny and Robson \cite{R1} (see also
\cite{R2}). With this in mind, it comes as no surprise that existence results
for subgame-perfect equilibria in continuous time are so hard to prove.

\section{Equilibrium strategies: definition\label{sec1}}

We consider an intertemporal decision problem where the decision-maker at time
$t$ is striving to maximise:
\begin{equation}
\int_{t}^{T}h\left(  s-t\right)  u\left(  c\left(  s\right)  \right)
ds+h\left(  T-t\right)  g\left(  k\left(  T\right)  \right)  \ \label{12}%
\end{equation}
subject to:%
\begin{equation}
\frac{dk}{ds}=f\left(  s,k\left(  s\right)  \right)  -c\left(  s\right)
,\ \ k\left(  t\right)  =k_{t}.\label{a13}%
\end{equation}

Notations are as stated in the introduction. Recall that $h$ is continuously
differentiable, with $h\left(  0\right)  =1,$ $h\left(  t\right)  \geq0$,
$h^{\prime}\left(  t\right)  \leq0$, and $h\left(  t\right)  \rightarrow0$
when $t\rightarrow\infty$. It will also be assumed that $u,f,g$ are twice
continuously differentiable, that $u$ is strictly concave, and that $f$ is
strictly concave with respect to $k$.

We shall denote by \thinspace$i:R^{d}\rightarrow R^{d}$ the inverse of the
derivative $u^{\prime}:R^{d}\rightarrow R^{d}:$%
\[
u^{\prime}\left(  c\right)  =x\Longleftrightarrow c=i\left(  x\right)
\]
and it will be assumed that it is continuously differentiable. We shall also
consider the Legendre-Fenchel transform $\tilde{u}$ of the concave function
$u$, defined by:%
\[
\tilde{u}(x)=\max_{c\in R^{d}}\left(  u(c)-xc\right)  =u\left(  i\left(
x\right)  \right)  -xi\left(  x\right)
\]

Note that it is a convex function. By the envelope theorem, we have:%
\[
\tilde{u}^{\prime}(x)=-i\left(  x\right)  =\left[  -u^{\prime}\right]
^{-1}\left(  x\right)
\]

We now proceed to define subgame-perfect equilibrium strategies, using the
approach outlined in the introduction. A strategy $c=\sigma\left(  t,k\right)
$ has been announced an is public knowledge. All decision-makers up to time
$t$ have applied this strategy, that is, the dynamics of capital between times
$0$ and $t$ are given by:%
\begin{equation}
\frac{dk}{ds}=f\left(  s,k\right)  -\sigma\left(  s,k\right)  ,\ \ k\left(
0\right)  =k_{0}\label{a1}%
\end{equation}

The decision-maker at time $t$ inherits a capital $k_{t}$, which is the value
at $s=t$ of the solution to the Cauchy problem (\ref{a1}). She can commit all
the decision-makers in $\left[  t,\ t+\varepsilon\right]  ,$where
$\varepsilon>0$ is vanishingly small. She expects all later ones to apply the
strategy $\sigma:[0,T]\times R^{d}\rightarrow R^{d}$, and she asks herself if
it is in her own interest to apply the same strategy, that is, to consume
$\sigma\left(  t,k\right)  $. If she consumes another bundle, $c$ say, the
immediate utility flow during $\left[  t,\ t+\varepsilon\right]  $ is
$u\left(  c\right)  \varepsilon$. At time $t+\varepsilon$, the resulting
capital will be $k+\left(  f\left(  t,k\right)  -c\right)  \varepsilon$, and
from then on, the strategy $\sigma$ will be applied. The consumption at time
$s\geq t+r$ is $c\left(  s\right)  =\sigma\left(  s,k\left(  s\right)
\right)  ,$ where
\begin{align}
&  \frac{dk}{ds}=f\left(  s,k\left(  s\right)  \right)  -\sigma\left(
s,k\left(  s\right)  \right)  ,\ \ s\geq t+\varepsilon\label{13}\\
&  k\left(  t+\varepsilon\right)  =k_{t}+\left(  f\left(  t,k_{t}\right)
-c\right)  \varepsilon\label{14}%
\end{align}

Denote by $k_{0}\left(  s\right)  $ the future path of capital if the
decision-maker at time $t$ applies the strategy $\sigma$, that is, if
$c=\sigma\left(  t,k_{t}\right)  $. The dynamic of $k_{0}$ is given by:
\begin{align}
\frac{dk_{0}}{ds}  &  =f\left(  s,k_{0}\left(  s\right)  \right)
-\sigma\left(  s,k_{0}\left(  s\right)  \right)  ,\ \ s\geq t\label{15}\\
k_{0}\left(  t\right)   &  =k_{t}\label{16}%
\end{align}

Write $k\left(  s\right)  =k_{0}\left(  s\right)  +k_{1}\left(  s\right)
\varepsilon$, plug that into (\ref{13}),(\ref{14}), keeping only terms of
first order in $\varepsilon$. We get:
\begin{align*}
\frac{dk_{0}}{ds}+\varepsilon\frac{dk_{1}}{ds}  &  =f\left(  s,k_{0}\left(
s\right)  \right)  +\varepsilon\frac{\partial f}{\partial k}\left(
s,k_{0}\left(  s\right)  \right)  k_{1}\left(  s\right)  -\sigma\left(
s,k_{0}\left(  s\right)  \right) \\
&  -\varepsilon\frac{\partial\sigma}{\partial k}\left(  s,k_{0}\left(
s\right)  \right)  k_{1}\left(  s\right)  ,\ \ s\geq t+\varepsilon,\\
k\left(  t+\varepsilon\right)   &  =k_{0}\left(  t+\varepsilon\right)
+\varepsilon k_{1}\left(  t+\varepsilon\right) \\
&  =k_{0}\left(  t\right)  +\varepsilon\frac{dk_{0}}{ds}\left(  t\right)
+\varepsilon k_{1}\left(  t+\varepsilon\right) \\
&  =k_{t}+\varepsilon\left(  f\left(  s,k_{t}\right)  -\sigma\left(
t,k_{t}\right)  \right)  +\varepsilon k_{1}\left(  t+\varepsilon\right)
\end{align*}
where $\frac{\partial f}{\partial k}$ and $\frac{\partial\sigma}{\partial k}$
stand for the matrix of partial derivatives of $f$ and $\sigma$ with respect
to $k\in R^{d}$. Comparing with (\ref{15}),(\ref{16}) and (\ref{14}), we get
the linear differential system:%
\begin{align*}
&  \frac{dk_{1}}{ds}=\left(  \frac{\partial f}{\partial k}\left(
s,k_{0}\left(  s\right)  \right)  -\frac{\partial\sigma}{\partial k}\left(
s,k_{0}\left(  s\right)  \right)  \right)  k_{1}\left(  s\right)  ,\ s\geq
t+\varepsilon\\
&  k_{1}\left(  t+\varepsilon\right)  =\sigma\left(  t,k_{t}\right)  -c
\end{align*}

Summing up, we find that the total gain for the decision-maker at time
$t\,\ $from consuming bundle $c$ during the interval of length $\varepsilon$
when she can commit, is:%
\begin{align*}
&  u\left(  c\right)  \varepsilon+\int_{t+\varepsilon}^{T}h\left(  s-t\right)
u\left(  \sigma\left(  s,k_{0}\left(  s\right)  +\varepsilon k_{1}\left(
s\right)  \right)  \right)  ds+h\left(  T-t\right)  \left(  g\left(
k_{0}\left(  T\right)  +\varepsilon k_{1}\left(  T\right)  \right)  \right) \\
&  =\int_{t}^{T}h\left(  s-t\right)  u\left(  \sigma\left(  s,k_{0}\left(
s\right)  \right)  \right)  ds+h\left(  T-t\right)  g\left(  k_{0}\left(
T\right)  \right) \\
&  +\varepsilon\left[  u\left(  c\right)  -u\left(  \sigma(t,k)\right)
+\int_{t}^{T}h\left(  s-t\right)  \frac{\partial u}{\partial c}\left(
\sigma\left(  s,k_{0}\left(  s\right)  \right)  \right)  \frac{\partial\sigma
}{\partial k}\left(  s,k_{0}\left(  s\right)  \right)  k_{1}\left(  s\right)
ds\right. \\
&  +\left.  h\left(  T-t\right)  \frac{\partial g}{\partial k}\left(
k_{0}\left(  T\right)  \right)  k_{1}\left(  T\right)  \right]  +\text{ h.o.t}%
\end{align*}
where $\frac{\partial g}{\partial k}$ (resp. $\frac{\partial u}{\partial c}$)
is the vector of partial derivatives of $g$ (resp. $u$) with respect to $k\in
R^{d}$ (resp. $c\in R^{d}$) and h.o.t denotes higher-order terms in
$\varepsilon$.

In the limit, when $\varepsilon\rightarrow0$, and the commitment span of the
decision-maker vanishes, we are left with two terms only. Note that the first
term does not depend on the decision taken at time $t$, but the second one
does. This is the one that the decision-maker at time $t$ will try to
maximize. In other words, given that a strategy $\sigma$ has been announced
and that the current state is $k_{t}=k\,$, the decision-maker at time
$t\,\ $faces the optimisation problem:%
\begin{equation}
\max_{c}P_{1}\left(  t,k,\sigma,c\right) \label{a14}%
\end{equation}
where:%
\begin{align*}
P_{1}\left(  t,k,\sigma,c\right)   &  =u\left(  c\right)  -u\left(
\sigma(t,k)\right) \\
&  +\int_{t}^{T}h\left(  s-t\right)  \frac{\partial u}{\partial c}\left(
\sigma\left(  s,k_{0}\left(  s\right)  \right)  \right)  \frac{\partial\sigma
}{\partial k}\left(  s,k_{0}\left(  s\right)  \right)  k_{1}\left(  s\right)
ds\\
&  +h\left(  T-t\right)  \frac{\partial g}{\partial k}\left(  k_{0}\left(
T\right)  \right)  k_{1}\left(  T\right)  ,
\end{align*}

In the above expression, $k_{0}\left(  s\right)  $ and $k_{1}\left(  s\right)
$ are given by:%
\begin{align*}
\frac{dk_{0}}{ds}  &  =f\left(  s,k_{0}\left(  s\right)  \right)
-\sigma\left(  s,k_{0}\left(  s\right)  \right) \\
k_{0}\left(  t\right)   &  =k
\end{align*}%
\begin{align}
\frac{dk_{1}}{ds}  &  =\left(  \frac{\partial f}{\partial k}\left(
s,k_{0}\left(  s\right)  \right)  -\frac{\partial\sigma}{\partial k}\left(
s,k_{0}\left(  s\right)  \right)  \right)  k_{1}\left(  s\right) \label{RI}\\
k_{1}\left(  t\right)   &  =\sigma\left(  t,k\right)  -c\label{RCI}%
\end{align}

\begin{definition}
We shall say that $\sigma:\left[  0,T\right]  \times R^{d}\rightarrow R^{d}$
is an \textbf{equilibrium strategy for the intertemporal decision model}
(\ref{12}),(\ref{a13}) if, for every $t\in\left[  0,T\right]  $ and $k\in
R^{d}$, the maximum in problem (\ref{a14}) is attained for $c=\sigma\left(
t,k\right)  $:%
\[
\sigma\left(  t,k\right)  =\arg\max_{c}P_{1}\left(  t,k,\sigma,c\right)
\]

\end{definition}

The intuition behind this definition is quite simple. Each decision-maker can
commit only for a small time $\varepsilon,\,$\ so he can only hope to exert a
very small influence on the final result. In fact, if the decision-maker at
time $t$ plays $c$ when he/she is called to bat, while all the others are
applying the strategy $\sigma$, the end payoff for him/her will be of the
form
\[
P_{0}\left(  t,k,\sigma\right)  +\varepsilon P_{1}\left(  t,k,\sigma,c\right)
\]
where the first term of the right hand side does not depend on $c$. In the
absence of commitment, the decision-maker at time $t$ will choose whichever
$c$ maximizes the second term $\varepsilon P_{1}\left(  t,k,\sigma,c\right)
$. Saying that $\sigma$ is an equilibrium strategy means that the decision
maker at time $t$ will choose $c=\sigma\left(  t,k\right)  $, that is, that
the strategy $\sigma$ can be implemented even in the absence of commitment.

Conversely, is a strategy $\sigma$ for the intertemporal decision model
(\ref{12}),(\ref{a13}) is not an equilibrium strategy, then it cannot be
implemented unless the decision-maker at time $0$ has some way to commit his
successors. Typically, an optimal strategy will not be an equilibrium
strategy. More precisely, a strategy which appears to be optimal at time $0$
no longer appears to be optimal at times $t>0$, which means that the
decision-maker at time $t$ feels he can do better than whatever was planned
for him to do at time $0$. In the case of macroeconomic policy, for instance,
successive governments will disagree on what is an optimal strategy, even if
they agree on the collective utility $u\left(  c\right)  $, so that the
concept of equilibrium strategy seems far more reasonable - at least it stands
a chance of being implemented.

What happens then if successive decision-makers take the myopic view, and each
of them acts as if he could commit his successors ?\ At time $t$, then, the
decision-maker would maximise the integral (\ref{12}) with the usual tools of
control theory, thereby deriving a consumption $c=\sigma_{n}\left(
t,k\right)  $. This is the \textit{naive strategy}; in general it will not be
an equilibrium strategy, so that every decison-maker has an incentive to
deviate. It will be studied in more detail in section \ref{sec3}.

\section{Characterization and existence of equilibrium strategies\label{sec2}}

In this section, we characterize equilibrium strategies of problem (\ref{10}),
(\ref{11}), (\ref{12}), by an equation, which we call the \emph{equilibrium
equation} (E), and which is reminescent - although different from - of the
Hamilton-Jacobi-Bellman (HJB) equation of optimal control. Note that there is
also an (HJB)\ equation associated with problem (\ref{3}),(\ref{5}), but it is
different from the equilibrium equation, and characterizes optimal strategies
instead of equilibrium ones. We will see that the only case when equations (E)
and (HJB) coincide is the case of exponential discount, and then equilibrium
strategies are also optimal strategies.

The equilibrium equation comes in two different guises:\ an integrated form
(IE) and a differentiated form (DE). We first derive the integrated form, and
then we show that it is equivalent to the differentiated one. Finally, under
suitable technical conditions on the utility function $u$ and the function
$f$, we show that solutions to the equilibrium equation exist close to the
terminal time $T$.

Given a strategy $\sigma\left(  t,k\right)  $, we shall be dealing with the
differential equation:%
\begin{align}
\frac{dk(s)}{ds}  &  =f\left(  s,k\left(  s\right)  \right)  -\sigma\left(
s,k\left(  s\right)  \right) \label{dynamics}\\
k\left(  t\right)   &  =k\nonumber
\end{align}

We shall denote by $\mathcal{K}\left(  s,t,k\right)  $ the flow associated
with this equation, that is the value at time $s$ of the solution of
(\ref{dynamics}) which takes the value $k$ at time $t$. It is defined by:
\begin{align}
\frac{\partial\mathcal{K}\left(  s,t,k\right)  }{\partial s}  &  =f\left(
s,\mathcal{K}\left(  s,t,k\right)  \right)  -\sigma\left(  s,\mathcal{K}%
\left(  s,t,k\right)  \right) \label{eq: HJBid}\\
\mathcal{K}\left(  t,t,k\right)   &  =k.\label{eq: HJBidic}%
\end{align}

In other words, $\mathcal{K}\left(  s,t,k\right)  $ is the value at time
$s\,$\ of the solution of:
\begin{equation}
\frac{dk}{ds}=f\left(  s,k\right)  -\sigma\left(  s,k\right) \label{25}%
\end{equation}
which takes the value $k$ at time $t$.

\subsection{Equilibrium characterization}

We shall say that a function $V:\left[  0,T\right]  \times R^{d}\rightarrow R$
satisfies \emph{the integrated equilibrium equation} (IE) if we have, for
every $t\in\left[  0,T\right]  $ and every $k:$%
\begin{equation}
V(t,k)=\int_{t}^{T}h(s-t)u\circ i\left(  \frac{\partial V}{\partial k}%
(s,k_{0}(s))\right)  ds+h(T-t)g(k_{0}(T))\tag{(IE)}%
\end{equation}
where:%
\begin{align*}
\frac{dk_{0}}{ds}  &  =f\left(  s,k_{0}\left(  s\right)  \right)  -i\circ
\frac{\partial V}{\partial k}(s,k_{0}(s))\\
k_{0}\left(  t\right)   &  =k
\end{align*}

Note that every solution of (IE)\ must satisfy the \emph{boundary condition}:%
\begin{equation}
V\left(  T,k\right)  =g\left(  k\right)  \ \ \ \forall k\tag{BC}%
\end{equation}

The following theorem characterizes the equilibrium strategies and its proof
is given in the Appendix A.

\begin{theorem}
\label{Th: condition necessaire}Let $\sigma:\left[  0,T\right]  \times
R^{d}\rightarrow R^{d}$ be jointly continuous, and continously differentiable
with respect to $k$ and let $\mathcal{K}$ be the associated flow defined by
(\ref{eq: HJBid}), (\ref{HJBidic}). Suppose $\sigma$ is an equilibrium
strategy for the intertemporal decision model (\ref{12}),(\ref{a13}). Then the
function:%
\begin{equation}
V(t,k)=\int_{t}^{T}h(s-t)u\left(  \sigma(s,\mathcal{K}\left(  s,t,k\right)
)\right)  ds+h(T-t)g\left(  \mathcal{K}\left(  T,t,k\right)  \right)
,\label{ct}%
\end{equation}
satisfies the integrated equilibrium equation (IE) and we have:%
\begin{equation}
\frac{\partial u}{\partial c}\left(  \sigma(t,k)\right)  =\frac{\partial
V}{\partial k}(t,k)\label{eq: Enveloppe}%
\end{equation}
Conversely, if a function $V$ is twice continuously differentiable and
satisfies the integrated equilibrium equation (IE), then:%
\[
\sigma(t,k)=i\left(  \frac{\partial V}{\partial k}(t,k)\right)
\]
is an equilibrium strategy.
\end{theorem}

Relation (\ref{eq: Enveloppe}) says that, along an equilibrium path, the
effect of an increment to current wealth on future utility, $\frac{\partial
V}{\partial k}(t,k)$, must balance the effect of an increment to current
consumption on current utility,$\frac{\partial u}{\partial c}\left(
\sigma(t,k)\right)  $. Thus, relation (\ref{eq: Enveloppe}) reflects the usual
tradeoff between utility derived from current consumption and utility value of saving.

From now on, we rewrite (IE) in the form
\begin{equation}
V(t,k)=\int_{t}^{T}h(s-t)u\left(  \sigma(s,\mathcal{K}\left(  s,t,k\right)
)\right)  ds+h(T-t)g\left(  \mathcal{K}\left(  T,t,k\right)  \right)  ,
\end{equation}
with the understanding that $\mathcal{K}\left(  s,t,k\right)  $ is the flow
associated with $\sigma\left(  t,k\right)  =i\circ\frac{\partial V}{\partial
k}\left(  t,k\right)  $.

The following proposition gives a differentiated version of the equilibrium equation.

\begin{proposition}
\label{Proposition: IE / DE} Assume that a function $V\left(  t,k\right)  $ is
twice continuously differentiable. Then $V$ satisfies the integrated
equilibrium equation (IE) if and only if it satisfies the differentiated
equilibrium equation:
\begin{align}
\frac{\partial V}{\partial t}(t,k)  &  +\int_{t}^{T}h^{\prime}(s-t)u\left(
\sigma(s,\mathcal{K}\left(  s,t,k\right)  )\right)  ds+h^{\prime
}(T-t)g(\mathcal{K}\left(  T,t,k\right)  )\nonumber\\
&  +\tilde{u}\left(  \frac{\partial V}{\partial k}(t,k)\right)  +\frac
{\partial V}{\partial k}(t,k)f(t,k)=0,\tag{DE}\label{dHJB}%
\end{align}
for all $(t,k)\in\lbrack0,T]\times R^{d}$, with the boundary condition
\begin{equation}
V(T,k)=g\left(  k\right)  .\tag{BC}\label{db}%
\end{equation}

\end{proposition}

Appendix B proves this proposition. It may be useful to rewrite it in the
following way:%
\begin{equation}
\rho(t,k)=\frac{1}{V}\left(  u\left(  \sigma(t,k)\right)  +\frac{\partial
V}{\partial t}(t,k)+\frac{\partial V}{\partial k}(t,k)\frac{\partial
\mathcal{K}}{\partial s}(t,t,k)\right) \label{dr}%
\end{equation}
where:%
\[
\rho(t,k)=-\frac{\int_{t}^{T}\frac{h^{\prime}(s-t)}{h(s-t)}h(s-t)u\left(
\sigma(s,\mathcal{K}\left(  s,t,k\right)  )\right)  ds+\frac{h^{\prime}%
(T-t)}{h(T-t)}h(T-t)g\left(  \mathcal{K}\left(  T,t,k\right)  \right)  }%
{\int_{t}^{T}h(s-t)u\left(  \sigma(s,\mathcal{K}\left(  s,t,k\right)
)\right)  ds+h(T-t)g\left(  \mathcal{K}\left(  T,t,k\right)  \right)  }%
\]
is interpreted as an effective discount rate. Equation (\ref{dr}) then tells
us that, along an equilibrium path, the relative changes in value to the
consumer must be equal to the effective discount rate.

Finally, when the discount rate is exponential, the effective discount rate is
just the constant discount rate $\rho=-h^{\prime}(t)/h(t)$ and equation (DE)
becomes simply the familiar (HJB) equation.

\begin{corollary}
\label{Proposition: HJB} With the exponential discounting $h(s)=e^{-\rho s}$,
the (DE) equation reduces to:
\begin{equation}
\frac{\partial V}{\partial t}(t,k)-\rho V(t,k)+\tilde{u}\left(  \frac{\partial
V}{\partial k}(t,k))\right)  +\frac{\partial V}{\partial k}%
(t,k)f(t,k)=0\label{eq: HJB}%
\end{equation}

\end{corollary}

\begin{proof}
In the exponential case, equations (IE)\ and (DE) become:%
\[
V(t,k)=\int_{t}^{T}e^{-\rho\left(  s-t\right)  }u\left(  \sigma(s,\mathcal{K}%
\left(  s,t,k\right)  )\right)  ds+e^{-\rho\left(  T-t\right)  }g\left(
\mathcal{K}\left(  T,t,k\right)  \right)  ,
\]%
\begin{align*}
\frac{\partial V}{\partial t}(t,k)  &  -\rho\int_{t}^{T}e^{-\rho\left(
s-t\right)  }u\left(  \sigma(s,\mathcal{K}\left(  s,t,k\right)  )\right)
ds-\rho e^{-\rho\left(  T-t\right)  }g(\mathcal{K}\left(  T,t,k\right)  )\\
&  +\tilde{u}\left(  \frac{\partial V}{\partial k}(t,k))\right)
+\frac{\partial V}{\partial k}(t,k)f(t,k)=0.
\end{align*}
Comparing, we immediately get (\ref{eq: HJB}).
\end{proof}

\subsection{Existence}

Neither equation (IE)\ nor equation (DE) are of a classical mathematical type.
If it were not for the integral term, equation (DE)\ would be a first-order
partial differential equation of known type (Hamilton-Jacobi), but this
additional term\ (an integral along the trajectory of the flow
(\ref{eq: HJBid}) associated with the solution $V\left(  t,k\right)  $ creates
additional complications.

In the sequel, we will solve that equation explicitly in particular cases. The
questions of existence and uniqueness in the general case are very much open.
In forhcoming work, Ekeland and Nirenberg prove a local existence result:

\begin{theorem}
Assume that all data ($u,f,g\,$\ and $h$) are analytic functions. Then, for
every $\bar{k}$, there are numbers $\varepsilon>0,\eta>0$ and a function
$V\left(  t,k\right)  $, defined for $T-\varepsilon\leq t\leq T$ and
$\left\Vert k-\bar{k}\right\Vert \leq\eta$, such that $V$ satisfies (DE)\ and (BC)
\end{theorem}

Recall that a function is analytic at a given point if its Taylor expansion at
that point has a non-zero radius of convergence. It is analytic if it is
analytic at every point. The proof of the theorem relies on a generalized
version of the classical Cauchy-Kowalewska theorem due to Nishida and
Nirenberg \cite{KN}.

\subsection{The infinite-horizon problem}

In the sequel, we will be looking at the infinite-horizon problem, whereby the
benefit to the decision-maker at time $t$ of a future consumption path
$s\rightarrow$ $c\left(  s\right)  ,s\geq t,$ is:%
\[
\int_{t}^{\infty}h\left(  s-t\right)  u\left(  c\left(  s\right)  \right)  ds
\]

The change of variables $s^{\prime}=s-t\geq0$ brings that integral to the
form:%
\[
\int_{0}^{\infty}h\left(  s^{\prime}\right)  u\left(  c\left(  s^{\prime
}+t\right)  \right)  ds^{\prime}%
\]
which is the benefit which the decision-maker at time $0$ derives from a
future consumption path $s^{\prime}\rightarrow c\left(  s^{\prime}+t\right)  $.

Assume now that the problem is stationary, meaning that the production
function $f\left(  t,k\right)  $ does not depend on $t:$%
\[
f\left(  t,k\right)  =f\left(  k\right)
\]

In that case, if the decision-maker at time $t$ resets his watch, so that time
$s$ becomes $s-t$, she faces exactly the same problem as the decision-maker at
time $0$. Under these circumstances, it is natural to expect that, if both
decision-makers have the same capital $k$, they will get the same equilibrium
value:%
\[
V\left(  t,k\right)  =V\left(  k\right)  \ \ \ \forall k
\]

We will now look directly for time-independent value functions. Consider the
equations:%
\begin{align}
V\left(  k\right)   &  =\int_{0}^{\infty}h\left(  t\right)  u\left(
i\circ\frac{\partial V}{\partial k}\left(  \mathcal{K}\left(  t,0,k\right)
\right)  \right)  dt\label{inf1}\\
0  &  =\frac{\partial V}{\partial k}(k)f(k)+\int_{0}^{\infty}h^{\prime
}(t)u\circ i\left(  \frac{\partial V}{\partial k}(\mathcal{K}\left(
t,0,k\right)  \right)  dt+\tilde{u}\left(  \frac{\partial V}{\partial
k}(k))\right) \label{inf2}%
\end{align}

\begin{lemma}
If a $C^{2}$ function $V\left(  k\right)  $ satisfies equation (\ref{inf1}) or
(\ref{inf2}), then $V\left(  t,k\right)  :=V\left(  k\right)  $ is a value
function for the infinite-horizon problem.
\end{lemma}

\begin{proof}
It is enough to show it for equation (\ref{inf1}). We have to prove that, for
every $t$, we have:
\[
V\left(  k\right)  =\int_{t}^{\infty}h\left(  s-t\right)  u\left(  i\circ
\frac{\partial V}{\partial k}\left(  \mathcal{K}\left(  s-t,t,k\right)
\right)  \right)  ds
\]
Note that the differential equation (\ref{eq: HJBid}) becomes autonomous, and
the function $\mathcal{K}$ displays the additional property
\[
\mathcal{K}\left(  s,t_{1},k\right)  =\mathcal{K}\left(  s-t_{2},t_{1}%
-t_{2},k\right)  ,~~~0<t_{2}<t_{1}<s.
\]
\

Changing variables in the integral, we get:%
\[
V\left(  k\right)  =\int_{0}^{\infty}h\left(  s\right)  u\left(  i\circ
\frac{\partial V}{\partial k}\left(  \mathcal{K}\left(  s,0,k\right)  \right)
\right)  ds
\]
which is precisely equation (\ref{inf1}).
\end{proof}

If $V\left(  k\right)  $ satisfies (\ref{inf1}) or (\ref{inf2}), the
corresponding equilibrium strategy:%
\[
\sigma\left(  k\right)  =i\circ\frac{\partial V}{\partial k}%
\]
which is time-independent, will be called \emph{stationary}. Note that a
stationary problem may have non-stationary equilibria.

\section{The consumption-saving problem\label{sec3}}

In this section, we assume $d=1$, so that there is only one good, and we will
be looking at a special case of the infinite-horizon problem. At any point in
time $s\in\lbrack0,\infty)$ the consumer has a stock of wealth $k(s)\in
(0,\infty)$ and receives a flow of labor income $w(s)$ as well as a flow of
interest income $r(s)k(s)$. Beginning with a capital stock $k\in(0,\infty)$ at
time $t$, we formulate the consumption-saving problem by
\begin{align}
&  \max\int_{t}^{\infty}h\left(  s-t\right)  u\left(  c\left(  s\right)
\right)  ds\ \label{Cs1}\\
&  \frac{dk(s)}{ds}=w(s)+r(s)k(s)-c(s),\ k\left(  t\right)  =k_{t}%
.\ \label{Cs3}%
\end{align}

This is a special case of the general problem (\ref{12}), (\ref{a13}), with:%
\[
f\left(  t,k\right)  :=w\left(  t\right)  +r\left(  t\right)  k
\]

We emphasize that at any point in time $t\in\lbrack0,\infty)$ the consumer
takes as given the interest rate $r(t)$ and the wage $w(t)$. Since $r\left(
t\right)  $ and $w\left(  t\right)  $ are time-dependent, we expect the value
function to be non-stationary, even though the horizon is infinite.

Equation (IE) becomes
\begin{equation}
V(t,k)=\int_{t}^{\infty}h(s-t)u\circ i\left(  \frac{\partial V}{\partial
k}(s,\mathcal{K}\left(  s,t,k\right)  )\right)  ds\label{IEcs}%
\end{equation}
where the flow $\mathcal{K}\left(  s,t,k\right)  $ solves
\begin{align}
\frac{d\mathcal{K}\left(  s,t,k\right)  }{ds}  &  =w(s)+r(s)\mathcal{K}\left(
s,t,k\right)  -i\circ\frac{\partial V}{\partial k}(s,\mathcal{K}\left(
s,t,k\right)  )\label{Fcs}\\
\mathcal{K}\left(  t,t,k\right)   &  =k.\label{IFcs}%
\end{align}

The next subsection gives explicit solutions when the utility function is in
the CRRA class.

\subsection{CRRA preferences}

In this section, we shall assume that $d=1$ and that the utility function
takes one of the forms:%
\begin{align*}
u(c)  &  =\frac{c^{1-\gamma}}{1-\gamma},~~\gamma>0\\
u\left(  c\right)   &  =\ln c
\end{align*}
\textit{\ }the latter corresponding to $\gamma=1$.An explicit construction of
the equilibrium strategy will be shown to be possible under an additional assumption:

\textbf{Assumption $A$:} \textit{There is at least one non-negative function
}$t\rightarrow\lambda\left(  t\right)  $\textit{\ which solves the fixed-point
problem:}%
\begin{equation}
\lambda(t)=\left[  \int_{t}^{\infty}\lambda(s)^{1-\gamma}\exp\left[  -\left(
1-\gamma\right)  \int_{t}^{s}(\lambda(u)-r(u))du\right]  h(s-t)ds\right]
^{-\frac{1}{\gamma}}\label{Recursion}%
\end{equation}

\begin{proposition}
\label{phe}If the utility function is CRRA and $\lambda\left(  t\right)  $ is
given by (\ref{Recursion}), the strategy
\begin{equation}
\sigma(t,k)=\lambda(t)\left[  k+\int_{t}^{\infty}\exp\left[  -\int_{t}%
^{s}r(u)du\right]  w(s)ds\right] \label{csep}%
\end{equation}
is an equilibrium strategy for the infinite-horizon problem. The associated
value function is given by
\[
V(t,k)=[\lambda(t)]^{-\gamma}\frac{\left[  k+\int_{t}^{\infty}\exp\left[
-\int_{t}^{s}r(u)du\right]  w(s)ds\right]  ^{1-\gamma}}{1-\gamma}.
\]

\end{proposition}

The equilibrium policy (\ref{csep}) consists of consuming the proportion
$\lambda\left(  t\right)  $ of current wealth; the latter is the sum of the
current capital stock and the present value of future wages.

The proof of Proposition \ref{phe} is given in Appendix \ref{G}

We now investigate equation (\ref{Recursion}) more closely. There are three
cases where it can be solved easily:

\subsubsection{Constant discounting}

In the case when $h(s)=e^{-\rho s}$, we find that the function
\[
\bar{\lambda}(t):=\frac{\exp\left[  \frac{1}{\gamma}\int_{0}^{t}\left(
(1-\gamma)r(s)-\rho\right)  ds\right]  }{\int_{t}^{\infty}\exp\left[  \frac
{1}{\gamma}\int_{0}^{s}\left(  (1-\gamma)r(u)-\rho\right)  du\right]  ds}%
\]
solves the equation (\ref{Recursion}) provided that the above integrals are
well-defined. Therefore, the policy
\begin{equation}
\sigma(t,k)=\bar{\lambda}(t)\left[  k+\int_{t}^{\infty}\exp\left[  -\int
_{t}^{s}r(u)du\right]  w(s)ds\right] \label{pcd}%
\end{equation}
is an equilibrium policy. Note that this is precisely the optimal policy from
the time $t$ perspective, which was expected anyway, since, with exponential
discount, optimal policies are equilibrium policies.

\subsubsection{Logarithmic utility}

We are now back with a general discount function $h\left(  t\right)  $, but we
choose a particular utility function, namely $u\left(  c\right)  =\ln c$, so
that $\gamma=1$. The equation reduces to:%
\[
\lambda(t)=\left[  \int_{t}^{\infty}h(s-t)ds\right]  ^{-1}=\left[  \int
_{0}^{\infty}h\left(  t\right)  dt\right]  ^{-1}%
\]

So $\lambda$ is constant\ (in spite of the fact that the interest rate on
capital $r\left(  t\right)  $ is time-dependent). This fact was first observed
by Barro \cite{Barro}. The corresponding equilibrium strategy is:
\begin{equation}
\sigma(t,k)=\frac{1}{\int_{0}^{\infty}h(s)ds}\left[  k+\int_{t}^{\infty}%
\exp\left[  -\int_{t}^{s}r(u)du\right]  w(s)ds\right]
\end{equation}

\subsubsection{Knife-edge case}

Assume the interest rate on capital is constant and given by
\[
r=\left(  \int_{0}^{\infty}h(s)ds\right)  ^{-1}%
\]

We seek to solve equation (\ref{Recursion}) for a constant $\lambda$. This yields:%

\begin{align*}
\lambda &  =\left[  \int_{t}^{\infty}\lambda^{1-\gamma}e^{-\left(
1-\gamma\right)  \left(  \lambda-r\right)  \left(  s-t\right)  }%
h(s-t)ds\right]  ^{-\frac{1}{_{\gamma}}}\\
&  =\left[  \int_{t}^{\infty}\lambda^{1-\gamma}h(s-t)ds\right]  ^{-\frac
{1}{_{\gamma}}}=\lambda^{-\frac{1-\gamma}{\gamma}}\left[  \int_{0}^{\infty
}h(s)ds\right]  ^{-\frac{1}{_{\gamma}}}%
\end{align*}
so that $\lambda=r$. The corresponding equilibrium is given by:
\[
\sigma(t,k)=r\left[  k+\int_{t}^{\infty}e^{-r(s-t)}w(s)ds\right]  .
\]

Given that along this equilibrium path, the consumers will consume the annuity
value of the wealth, the above equilibrium is consistent with Friedman's
permanent-income model. Note that this equilibrium strategy would also be the
optimal strategy for the case of a constant discount rate equal to $r,$ so
that $h(s)=e^{-rs}$.

\subsection{Constant interest rate}

In this subsection, we shall assume that the interest rate on capital is
constant:
\[
r\left(  t\right)  =r
\]

The following examples provide, for some specific discount functions $h$,
explicit formulas for some equilibrium strategies with constant propensity to
consume out of wealth.

\subsubsection{Exponential discount}

When the discount function is exponential $h(s)=e^{-\rho s}$, the equation
(\ref{Recursion}) for $\lambda$ takes the form
\[
1=\lambda\int_{t}^{\infty}e^{-\left(  \rho+(\lambda-r)(1-\gamma)\right)
(s-t)}ds,
\]
where $r>0$ is the constant interest rate. If the above integral is well
defined, we see that:
\[
\lambda=r+\frac{\rho-r}{\gamma}=:\lambda_{0}.
\]

The policy
\[
\sigma(t,k)=\lambda_{0}\left[  k+\int_{t}^{\infty}e^{-r\left(  s-t\right)
}w(s)ds\right]
\]
is an equilibrium policy provided $\lambda_{0}>0$, that is $\rho
-r(1-\gamma)>0$.

\subsubsection{A mixture of exponential discount functions}

For convenience, assume that $\gamma>1$.

Consider the case when the discount function is the mixture of two exponential
functions, that is:
\[
h(s)=\omega e^{-\rho_{1}s}+(1-\omega)e^{-\rho_{2}s},
\]
where $0<\rho_{1}<\rho_{2}$, and $\omega\in\lbrack0,1]$. The instantaneous
discount rate associated to $h$ at time $t$ is
\[
-\frac{h^{\prime}(t)}{h(t)}=\rho_{1}+(\rho_{2}-\rho_{1})\frac{1-\omega}{\omega
e^{(\rho_{2}-\rho_{1})t}+(1-\omega)}%
\]
and is gradually declining from $\rho_{0}:=\omega\rho_{1}+(1-\omega)\rho_{2}$
(at time $t=0$) to $\rho_{1}$ (at time $t=\infty$). Therefore, this
specification captures the idea that discount rates decline with the horizon
over which utility is discounted, a feature that O'Donoghue an Rabin
(\cite{O'Donoghue-Rabin1}, \cite{O'Donoghue-Rabin2}) call the
\textquotedblleft present bias". The mixture of exponential discount function
also corresponds to what Harris and Laibson \cite{Harris-Laibson2} call the
"auxiliary model".

If the discount rate were constant and equal to the long term value $\rho_{1}
$, we would have a stationary equilibrium policy where ${\underline{\lambda}%
}=r+\left(  \rho_{1}-r\right)  /\gamma.$If the discount rate were constant and
equal to the short term value $\rho_{0}$, we would have a stationary
equilibrium policy where ${\overline{\lambda}}=r+\left(  \rho_{0}-r\right)
/\gamma$. Each of them would be optimal in its own context, given that the
discount rate is constant.

In the general case where $\omega\in(0,1)$, so that the discount rate declines
from $\rho_{0}$ to $\rho_{1}$, we look for an equilibrium policy where the
propensity to consume out of wealth is a constant $\lambda$. After
integrating, equation (\ref{Recursion}) turns out to be equivalent to the
following:
\begin{equation}
f(\lambda):=\frac{\omega}{\rho_{1}+(\lambda-r)(1-\gamma)}+\frac{1-\omega}%
{\rho_{2}+(\lambda-r)(1-\gamma)}-\frac{1}{\lambda}%
=0,\label{eq: mixture recursion}%
\end{equation}
provided that the integrability conditions:
\begin{equation}
\rho_{i}+(\lambda-r)(1-\gamma)>0,~~~i=1,2\label{eq: mixture integrability}%
\end{equation}
are satisfied.

The function $f$ is increasing on the interval $(0,r+\frac{\rho_{1}}{\gamma
-1})$ and furthermore $f(0)=-\infty$ and $f(r+\frac{\rho_{1}}{\gamma
-1})=+\infty$. Therefore there must exist a unique value $\lambda_{1}%
\in(0,r+\frac{\rho_{1}}{\gamma-1})$ such that $f(\lambda_{1})=0$. Recalling
that $\gamma>1$, it is easy to see that $\lambda_{1}$ satisfies the
integrability condition (\ref{eq: mixture integrability}) and a further
inspection reveals that $\lambda_{1}$ is the unique solution of the recursion
(\ref{eq: mixture recursion}) satisfying the integrability condition
(\ref{eq: mixture integrability}). Therefore, $\lambda_{1}$ gives rise to an
equilibrium strategy.

Evaluating $f$ at ${\underline{\lambda}}$ gives
\[
f({\underline{\lambda}})=\gamma(1-\gamma)\left[  \frac{1}{\rho_{1}+\gamma
(\rho_{2}-\rho_{1})-r(1-\gamma)}-\frac{1}{\rho_{1}-r(1-\gamma)}\right]  <0,
\]
and since $f$ is increasing, we obtain that ${\underline{\lambda}}<\lambda
_{1}$.

If the interest rate $r$ has the precise value:
\[
r=\frac{1}{\frac{\omega}{\rho_{1}}+\frac{1-\omega}{\rho_{2}}}%
\]
then $\lambda_{1}=r$ is the solution, and in that case ${\underline{\lambda}%
}<\lambda_{1}<$ $\bar{\lambda}.$

\subsubsection{Quasi hyperbolic discount.}

We define the discount function (in continuous time) as
\[
h(s)=\left\{
\begin{array}
[c]{ll}%
e^{-\rho s} & \text{for }0\leq s\leq\tau\\
\delta e^{-\rho s} & \text{for }s>\tau
\end{array}
\right.
\]
where $\tau>0$ and $\delta\in(0,1]$.

When $\delta=1$, the discount is exponential and the equilibrium is the one
described in the preceding subsection assuming that $\tilde{\rho}%
=\rho-r(1-\gamma)>0$. If $\delta<1$, assuming a time invariant propensity to
consume out of wealth, and integrating the equation (\ref{Recursion}) yields:
\[
1=\frac{\lambda}{\tilde{\rho}+\lambda(1-\gamma)}\left[  1-(1-\delta
)e^{-\left(  {\tilde{\rho}}+\lambda(1-\gamma)\right)  \tau}\right]  ,
\]
or equivalently $f\left(  \lambda\right)  =0$, where:%
\[
f\left(  \lambda\right)  =\gamma-(1-\delta)e^{-\left(  {\tilde{\rho}}%
+\lambda(1-\gamma)\right)  \tau}-\frac{\tilde{\rho}}{\lambda}%
\]
provided that an integrability condition holds:
\begin{equation}
\tilde{\rho}+\lambda(1-\gamma)>0\label{eq: Integrability}%
\end{equation}

If $0<\gamma<1$, we see that $f^{\prime}>0$ and thus $f$ is non decreasing. On
the other hand, when $f(\lambda)\rightarrow-\infty$ when $\lambda
\rightarrow0,\lambda\geq0$, and $f(\lambda)\rightarrow\gamma$ when
$\lambda\rightarrow\infty\,\ $so that there must exist a unique $\lambda_{2}$
such that $f(\lambda_{2})=0$. Since $\gamma<1$, the integrability condition
(\ref{eq: Integrability}) is satisfied

If $\gamma>1,$ we have $f\left(  \lambda\right)  \rightarrow-\infty\,$\ when
$\lambda\rightarrow0,\lambda>0$, and $f\left(  -\frac{\tilde{\rho}}{1-\gamma
}\right)  =\delta>0,$ so that $f$ has at least one root $\lambda_{2}$
satisfying the integrability condition.

So the existence of an equilibrium strategy with a constant propensity to
consume $\lambda\left(  t\right)  =\lambda_{2}$ is proved in all cases.

In the limiting case when $\tau\rightarrow0$, whe obtain the instant
gratification model of Harris and Laibson \cite{Harris-Laibson2}. Then:
\[
\lambda_{2}\rightarrow\frac{\rho+r(\gamma-1)}{\delta+\gamma-1}%
\]
which satisfies the integrability condition (\ref{eq: Integrability}) when
$\gamma>1$.

\subsubsection{General hyperbolic discount function}

We consider the discount function
\[
h(s)=\frac{1}{\left(  1+\alpha s\right)  ^{\frac{\beta}{\alpha}}}e^{-\rho
s},~~\alpha>0,\beta>0\text{ \ and }\rho>0
\]
specified by Luttmer and Mariotti \cite{Luttmer-Mariotti} and which
particularizes, when $\rho=0$, the generalized hyperbolic discount function
reported in Loewenstein and Prelec \cite{Loewenstein-Prelec}. The resulting
discount rate:
\[
-\frac{h^{\prime}(s)}{h(s)}=\rho+\frac{\beta}{1+\alpha s}%
\]
is smoothly declining from $\rho+\beta$ (at $s=0$) to $\rho$ (at $s=\infty$).
The coefficient $\alpha$ determines how close the discount function $h$ is to
the exponentials $e^{-\rho s}$ and $e^{-(\rho+\beta)s}$.

The equation (\ref{Recursion}) becomes
\[
1=\lambda\int_{0}^{\infty}\frac{1}{\left(  1+\alpha s\right)  ^{\frac{\beta
}{\alpha}}}e^{-(\rho+(\lambda-r)(1-\gamma))s}ds
\]
provided the integrability condition $\rho+(\lambda-r)(1-\gamma)>0$ is satisfied.

We define the function
\[
f(\lambda)=\int_{0}^{\infty}\frac{1}{\left(  1+\alpha s\right)  ^{\frac{\beta
}{\alpha}}}e^{-(\rho+(\lambda-r)(1-\gamma))s}ds-\frac{1}{\lambda}%
\]
and verify that $f(0)=-\infty$, $f\left(  r+\frac{\rho}{\gamma-1}\right)
=+\infty$ and $f^{\prime}>0$. Therefore, there exist a unique
\[
\lambda_{3}\in(0,r+\frac{\rho}{\gamma-1})
\]
such that $f(\lambda_{3})=0$ and such that the integrability condition
$\rho+(\lambda_{3}-r)(1-\gamma)>0$ is satisfied.

\subsection{Comparative analysis.}

We want compare the equilibrium strategy with the strategy which, from the
point of view of the decision-maker at time $t=0$, is optimal. We shall do so
in the case when the interest rate $r$ and the wage $w$ are constant, and when
$u\left(  c\right)  =\ln c.$

The equilibrium strategy, as we saw earlier, then is time-independent and
consists of consuming a constant fraction of current wealth:%
\begin{equation}
\sigma(k)=\frac{1}{\int_{0}^{\infty}h(s)ds}\left[  k+\frac{w}{r}\right]
\label{ek1}%
\end{equation}

Note that, for the model to be meaningful, we must have:%
\[
r>\frac{1}{\int_{0}^{\infty}h(s)ds}%
\]
otherwise equation (\ref{ek1}) would mean that in equilibrium, consumption is
greater that income. This makes sense: if the interest on capital is lower
than the psychological discount rate, there is no point in investing.

Let us put ourselves in the shoes of the decision-maker at time $t=0$, endowed
with a capital $k_{0}$, and find the optimal strategy from her point of view.
Solving the optimal control problem:%
\begin{align*}
&  \max\int_{0}^{\infty}h\left(  t\right)  \ln c\left(  t\right)  dt\\
\frac{dk}{dt}  &  =rk+w-c,\ \ k\left(  0\right)  =k_{0}%
\end{align*}
we find, by the Euler-Lagrange equation:%
\[
\frac{1}{c}\frac{dc}{dt}=r+\frac{h^{\prime}\left(  t\right)  }{h\left(
t\right)  }%
\]
which we integrate, to get:%
\[
c\left(  t\right)  =c_{0}h\left(  t\right)  e^{rt}%
\]

Substituting into the dynamics, we get:%
\[
\frac{dk}{dt}=rk+w-c_{0}h\left(  t\right)  e^{rt}%
\]
which we integrate, to get:%
\[
k\left(  t\right)  =\left(  k_{0}+\frac{w}{r}-c_{0}\int_{0}^{t}h\left(
s\right)  ds\right)  e^{rt}-\frac{w}{r}%
\]

Because of the transversality condition at infinity, we must have:%
\begin{equation}
c_{0}=\frac{1}{\int_{0}^{\infty}h\left(  s\right)  ds}\left(  k_{0}+\frac
{w}{r}\right) \label{ek3}%
\end{equation}
and the optimal propensity to consume at time $t$ is:%
\begin{align*}
\lambda\left(  t\right)   &  =\frac{c\left(  t\right)  }{k\left(  t\right)
+\frac{w}{r}}=\frac{c_{0}h\left(  t\right)  }{\left(  k_{0}+\frac{w}{r}%
-c_{0}\int_{0}^{t}h\left(  s\right)  ds\right)  }\\
&  =\frac{h\left(  t\right)  }{\left(  H-\int_{0}^{t}h\left(  s\right)
ds\right)  }=h\left(  t\right)  \left(  \int_{t}^{\infty}h\left(  s\right)
ds\right)  ^{-1}%
\end{align*}

At time $t=0$, we find $\lambda\left(  0\right)  =\left(  \int_{0}^{\infty
}h\left(  s\right)  ds\right)  ^{-1}$. This is precisely the equilibrium
value, as defined by (\ref{ek1}). The optimal propensity to consume,
$\lambda\left(  t\right)  $, is time-dependent, and deviates from the
equilibrium value $\lambda\left(  0\right)  $. Note that:%
\[
\lambda^{\prime}\left(  0\right)  =\frac{1}{\int_{0}^{\infty}h\left(
s\right)  ds}\left(  h^{\prime}\left(  0\right)  +\frac{1}{\int_{0}^{\infty
}h\left(  s\right)  ds}\right)
\]
so that $\lambda\left(  t\right)  $ may be greater or smaller than the
equilibrium value, according to the characteristics of the discount function
$h$. Applying the optimal strategy (from the point of view of time $0$) yields
the following dynamic (we set $\int_{0}^{\infty}h\left(  s\right)  ds=H$ for
the sake of convenience) ::%
\begin{align*}
k\left(  t\right)   &  =\left(  k_{0}+\frac{w}{r}-\left(  k_{0}+\frac{w}%
{r}\right)  \frac{1}{H}\int_{0}^{t}h\left(  s\right)  ds\right)  e^{rt}%
-\frac{w}{r}\\
c\left(  t\right)   &  =c_{0}h\left(  t\right)  e^{rt}=\frac{1}{H}\left(
k_{0}+\frac{w}{r}\right)  h\left(  t\right)  e^{rt}%
\end{align*}

Applying the equilibrium strategy yields the following dynamics:
\begin{align*}
k\left(  t\right)   &  =k_{0}e^{\left(  r-1/H\right)  t}+w\frac{1-1/rH}%
{r-1/H}\left(  1-e^{\left(  r-1/H\right)  t}\right) \\
c\left(  t\right)   &  =\frac{1}{H}\left(  \frac{w}{r}+k_{0}e^{\left(
r-1/H\right)  t}+w\frac{1-1/rH}{r-1/H}\left(  1-e^{\left(  r-1/H\right)
t}\right)  \right)
\end{align*}

Note, however, a remarkable fact. Define the \textit{naive strategy} as
follows: every decision-maker acts as if she could commit her successors; she
computes the control \ $c\left(  s\right)  $ which is optimal on the interval
$\left[  t,\infty\right]  $, and consumes $c\left(  t\right)  $. From the
previous analysis it follows that the naive strategy is an equilibrium
strategy. This, of course, is particular to the logarithmic case $u\left(
c\right)  =\ln c$.

\section{Indeterminacy in the Ramsey growth problem\label{sec4}}

We now go back to the general problem (\ref{12}), (\ref{a13}) in the
stationary case, where the production function is given by:%
\[
f\left(  t,k\right)  =f\left(  k\right)
\]

We then interpret the problem as the Ramsey problem in growth theory. It is
well-known, and described for instance in the textbook by Barro and
Sala-i-Martin \cite{Barro salai}, that there are two versions to that problem:

\begin{enumerate}
\item \emph{The centralized version}. A benevolent planner, seeking to
maximize the integral (\ref{12}) which measures global welfare, determines an
optimal growth strategy $\sigma\left(  t,k\right)  $ and commits citizens to
consume $\sigma\left(  t,k\right)  $ and to invest the remainder

\item \emph{The decentralized version}. Future interest rates $r_{1}\left(
t\right)  $ and future wages $w_{1}\left(  t\right)  $ are common knowledge.
The representative individual then solves the consumption-saving problem,
taking $r_{1}\left(  t\right)  $ and $w_{1}\left(  t\right)  $ as given. This
determines the rate of investment at any time $t$. This is turn determines the
wages $w_{2}\left(  t\right)  $ that the production sector can offer, and the
capital rental interest $r_{2}\left(  t\right)  $. One wants $r_{1}=r_{2}$ and
$w_{1}=w_{2}$.
\end{enumerate}

In the case of exponential discount, $h\left(  t\right)  =e^{-\rho t}$, both
problems have the same solution (see \cite{Barro salai}): if $\sigma\left(
t,k\right)  $ is a solution of the centralized problem, and $k\left(
t\right)  $ the corresponding optimal trajectory, then $w_{1}\left(  t\right)
:=f\left(  k\right)  -k\left(  t\right)  f^{\prime}\left(  k\left(  t\right)
\right)  $ and $r_{1}\left(  t\right)  :=f^{\prime}\left(  k\right)  $ have
the property that $r_{1}=r_{2}$ and $w_{1}=w_{2}$.

These notions naturally extend to more general discount functions. In the
absence of commitment technology, one must replace optimal policies by
equilibrium policies, and one is naturally led to two notions of and
equilibrium growth policy:

\begin{enumerate}
\item \emph{The centralized version}. There is a succession of benevolent
planners, each of them holding power during an infinitesimal period of time,
and having the ability to commit their contemporaries\ in the consumption and
production sectors during that period. They agree on an equilibrium strategy
$\sigma\left(  t,k\right)  $ for the problem (\ref{12}), (\ref{a13}).

\item \emph{The decentralized version}. Future interest rates $r_{1}\left(
t\right)  $ and future wages $w_{1}\left(  t\right)  $ are common knowledge.
There is a succession of representative individuals, and they agree on an
equilibrium strategy for the consumption-saving problem (\ref{Cs1}),
(\ref{Cs3}), taking $r_{1}\left(  t\right)  $ and $w_{1}\left(  t\right)  $ as
given. This determines the rate of investment at any time $t$. This is turn
determines the wages $w_{2}\left(  t\right)  $ that the production sector can
offer, and the capital rental interest $r_{2}\left(  t\right)  $. One wants
$r_{1}=r_{2}$ and $w_{1}=w_{2}$.
\end{enumerate}

In continuity with the results for the exponential discount, one would
naturally expect that the two problems coincide, but this is no longer the case.

The results in Barro \cite{Barro} pertain to the second problem. To the best
of our knowledge the first one, that is, the study of the planner's problem in
optimal growth theory under time inconsistency, has not been studied. The
remainder of this paper is devoted to shedding some light on that problem. As
in classical growth theory, we will concentrate on the one-dimensional
case:\ $d=1$ .

\begin{definition}
Take a point $\bar{k}$. We shall say that $\bar{k}$ is an \emph{equilibrium
point} if there is a stationary equilibrium strategy $\sigma\left(  k\right)
$, defined on a neighbourhood $\ \Omega$ of $\bar{k}$ in $R^{d}$, and such
that all trajectories of (\ref{dynamics}) starting inside $\Omega$ when $t=0$
converge to $\bar{k}$ when $t\rightarrow\infty$.
\end{definition}

It follows from the definition that the trajectory starting from $\bar{k}$ is
$\bar{k}$ itself:\ the solution of (\ref{dynamics}) with $k\left(  0\right)
=\bar{k}$ is $k\left(  t\right)  =\bar{k}$ for all $t$. Denoting by $\bar{c}$
the consumption along that trajectory, we must have:%
\[
\bar{c}=f\left(  \bar{k}\right)
\]

\begin{theorem}
\label{Theorem: indeterminacy} Assume that $\bar{k}$ is an equilibrium point,
and that the corresponding value function $V\left(  k\right)  $ is $C^{2}$ in
a neighbourhood of $\bar{k}$. Then the number $\alpha$ defined by:
\[
\alpha:=\frac{V^{\prime\prime}\left(  \bar{k}\right)  }{u^{\prime\prime
}\left(  \bar{c}\right)  }%
\]
must satisfy:%
\begin{equation}
\alpha\geq f^{\prime}\left(  \bar{k}\right)  >0\label{59}%
\end{equation}

If $\alpha>f^{\prime}\left(  \bar{k}\right)  $, then:%
\begin{equation}
\alpha\int_{0}^{\infty}h\left(  t\right)  \exp\left[  \left(  f^{\prime
}\left(  \bar{k}\right)  -\alpha\right)  t\right]  dt=1\label{51}%
\end{equation}

\end{theorem}

The proof is given in Appendix \ref{H}

\begin{corollary}
\label{corrolary: exponential discount} Set $h\left(  t\right)  =e^{-\rho t}$.
Assume that $\bar{k}$ is an equilibrium point. Then:%
\begin{equation}
f^{\prime}\left(  \bar{k}\right)  =\rho\label{53}%
\end{equation}
and,
\begin{equation}
\alpha\equiv\frac{V^{\prime\prime}(\bar{k})}{u^{\prime\prime}(\bar{c}%
)}=f^{\prime}(\bar{k})\left[  \frac{1+\sqrt{1+4\frac{u^{\prime}(\bar
{k})f^{\prime\prime}(\bar{k})}{\rho^{2}u^{\prime\prime}(\bar{k})}}}{2}\right]
>f^{\prime}(\bar{k})\label{53bis}%
\end{equation}

\end{corollary}

In the exponential case equation (\ref{51}) degenerates: it sets no condition
on $\alpha$, but determines $\bar{k}$ through (\ref{53}). This is the
well-known relation for the optimal growth path, which usually is obtained by
the transversality condition at infinity, and which here is derived in a novel way.

In the general case, as we will see in the following example, equation
(\ref{51}) does not determine $\bar{k}$: it determines $\alpha$ as a function
of $\bar{k}$. The proof is given in the appendix.

\begin{proposition}
\label{pol}Set $h\left(  t\right)  =e^{-rt}$ on $\left[  0,\ T\right]  $, and
$h\left(  t\right)  =0$ for $t>T$. Assume that $\bar{k}$ is an equilibrium
point. Then, there exists a decreasing function $\varphi:]0,\infty
\lbrack\rightarrow]0,\infty\lbrack$, with:%
\begin{align*}
\varphi\left(  \alpha\right)   &  \rightarrow\infty\text{ when }%
\alpha\rightarrow0\\
\varphi\left(  1\right)   &  =1
\end{align*}
and a number $a\left(  T\right)  \leq1/T$ such that conditions (\ref{51}),
(\ref{59}) are equivalent to the following:%
\begin{align}
0  &  <f^{\prime}\left(  \bar{k}\right)  -r\leq1/T\label{55}\\
\alpha T  &  =\varphi\left(  \left[  f^{\prime}\left(  \bar{k}\right)
-r\right]  T\right)  \text{ }\label{62}%
\end{align}

\end{proposition}

In other words, there is a continuous family of solutions to equation
(\ref{51}), one for each $\bar{k}$ such that $f^{\prime}\left(  \bar
{k}\right)  $ falls in the interval $(\rho,\ \rho+1/T)$. The corresponding
$\alpha\left(  \bar{k}\right)  $ goes to $\infty$ when $f^{\prime}\left(
\bar{k}\right)  \rightarrow\rho$ and to $1/T$ when $f^{\prime}\left(  \bar
{k}\right)  \rightarrow\rho+1/T$. Note that there is no solution for
$f^{\prime}\left(  \bar{k}\right)  \leq\rho$ of $f^{\prime}\left(  \bar
{k}\right)  \geq\rho+1/T$

For future reference, we write a few properties of the function $\varphi$.
They follow easily from the properties of the function $xe^{-x}:$%
\begin{align*}
x_{1}  &  <x_{2}\Longrightarrow\varphi\left(  x_{1}\right)  <\varphi\left(
x_{2}\right) \\
\varphi\left(  x\right)   &  \rightarrow\infty\text{ when }x\rightarrow
\infty,\\
\varphi\left(  x\right)   &  \rightarrow0\text{ when }x\rightarrow\infty\\
x\varphi\left(  x\right)   &  =\psi\left(  xe^{-x}\right)  \text{ with }%
\psi\left(  0\right)  =0\text{, }\psi^{\prime}\left(  0\right)  =1,\text{when
}x\geq1
\end{align*}

We now want to know what happens when the horizon $T$ goes to $0$ (instant
gratification) or $\infty$ (exponential discount). Let $\left(  \alpha\left(
T\right)  ,k\left(  T\right)  \right)  $ be a solution of equation (\ref{51})
with $h\left(  t\right)  $ as above:%
\begin{equation}
\alpha\left(  T\right)  \int_{0}^{T}e^{-\rho t}\exp\left[  \left(  f^{\prime
}\left(  k\left(  T\right)  \right)  -\alpha\left(  T\right)  \right)
t\right]  dt=1
\end{equation}

\begin{itemize}
\item Let $T\rightarrow\infty$.\ Assume $k\left(  T\right)  \rightarrow\bar
{k}$ and $\alpha\left(  T\right)  \rightarrow\bar{\alpha}>0$. Then:%
\[
f^{\prime}\left(  \bar{k}\right)  =\rho
\]

\item Let $T\rightarrow0$.\ Assume $k\left(  T\right)  \rightarrow\bar{k}$ and
$\alpha\left(  T\right)  \rightarrow\bar{\alpha}>0$. Then:%
\[
\bar{k}=0
\]

\end{itemize}

The proof of the first part follows immediately from the estimate
$0<f^{\prime}\left(  k\left(  T\right)  \right)  -\rho\leq1/T$. For the
second, we have $\alpha\left(  T\right)  =\varphi\left(  f^{\prime}\left(
k\left(  T\right)  \right)  T-\rho T\right)  /T$. Since the left-hand side
converges, so must the right-hand side, so $\varphi\left(  f^{\prime}\left(
k\left(  T\right)  \right)  T-\rho T\right)  $ must go to infinity when
$T\rightarrow0$, which is only possible if $f^{\prime}\left(  k\left(
T\right)  \right)  \rightarrow\infty$. Since $f$ satisfies the Inada
conditions, we must have $k\left(  T\right)  \rightarrow0,$ as announced,

These results are conform to economic intuition. Note in particular that, when
$T\rightarrow\infty$, we find again the condition $f^{\prime}\left(  \bar
{k}\right)  =\rho$ in the limit. However, for finite $T$, equation (\ref{51})
does not determine $k\left(  T\right)  $, which is a striking difference with
$T=\infty$. In this, as in the general case of non-constant discount, we have
been unable to find any further condition that would determine $k\left(
T\right)  $. This would indicate non-uniqueness of possible $k\left(
T\right)  $, and hence a multiplicity of equilibrium strategies, one for each
possible value of $k\left(  T\right)  $ and $\alpha\left(  T\right)  $. The
following results, which are valid in the case of general non-constant
discounts, indicates that, even with non-uniqueness, there is a definite range
of possible equilibrium values for $\bar{k}$.

\begin{corollary}
\label{corrolary: robustess 1} Assume that the discount function $h\left(
t\right)  $ satisfies:%
\[
e^{-\rho_{2}t}\leq h\left(  t\right)  \leq e^{-\rho_{1}t}\text{ for }t\geq0
\]
for some $0<\rho_{1}\leq\rho_{2}$, and that $\bar{k}$ is an equilibrium point
for $h$. Then:%
\[
\rho_{1}\,\leq f^{\prime}\left(  \bar{k}\right)  \leq\rho_{2}%
\]

\end{corollary}

\begin{proof}
From equation (\ref{51}), we have:%
\begin{align}
1  &  =\alpha\int_{0}^{\infty}h\left(  t\right)  \exp\left[  \left(
f^{\prime}\left(  \bar{k}\right)  -\alpha\right)  t\right]  dt\\
&  \leq\alpha\int_{0}^{\infty}\exp\left[  \left(  -r_{1}+f^{\prime}\left(
\bar{k}\right)  -\alpha\right)  t\right]  dt\\
&  =\frac{\alpha}{r_{1}-f^{\prime}\left(  \bar{k}\right)  +\alpha}%
\end{align}
and this gives $f^{\prime}\left(  \bar{k}\right)  \geq r_{1}$. Hence the result
\end{proof}

More generally, we have the following result:

\begin{proposition}
\label{proposition: robustess 2} Assume that the discount functions $h_{1}$
and $h$ satisfy:%
\[
h_{0}\left(  t\right)  \leq h_{1}\left(  t\right)  \text{ \ for }t\geq0
\]
Denote by $K_{0}$ and~$K_{1}$ the set of equilibrium points for $h_{0}$ and
$h_{1}$ respectively. Then:
\[
\sup K_{0}\leq\sup K_{1}%
\]

\end{proposition}

\begin{proof}
Assume otherwise, so that there exists some $\bar{k}_{0}\in K_{0}$ such that
$\bar{k}_{0}\,>\sup K_{1}$. Then there exists some $\bar{k}_{1}\notin K_{1}$
with $\bar{k}_{1}<\bar{k}_{0}$. . Since $f$ is strictly concave, we must have
$f^{\prime}\left(  \bar{k}_{0}\right)  <f^{\prime}\left(  \bar{k}_{1}\right)
$. Set:%
\begin{align}
\varphi_{0}\left(  \alpha\right)   &  :=\alpha\int_{0}^{\infty}h_{0}\left(
t\right)  \exp\left[  \left(  f^{\prime}\left(  \bar{k}_{0}\right)
-\alpha\right)  t\right]  dt\label{60}\\
\varphi_{1}\left(  \alpha\right)   &  :=\alpha\int_{0}^{\infty}h_{1}\left(
t\right)  \exp\left[  \left(  f^{\prime}\left(  \bar{k}_{1}\right)
-\alpha\right)  t\right]  dt\label{61}%
\end{align}

Since $\bar{k}_{1}\notin K_{1}$, we have $\varphi_{1}\left(  \alpha\right)
\neq1$ for all $\alpha>0$. Since $\varphi_{1}\left(  \alpha\right)
\rightarrow0$ when $\alpha\rightarrow\infty$, this implies that $\varphi
_{1}\left(  \alpha\right)  <1$ for all $\alpha>0$. Since $h_{0}\leq h_{1}$ and
$f^{\prime}\left(  \bar{k}_{0}\right)  \leq f^{\prime}\left(  \bar{k}%
_{1}\right)  $, we have $\varphi_{0}\leq\varphi_{1}$, such that $\varphi
_{1}\left(  \alpha\right)  <1$ for all $\alpha>0$.
\end{proof}

\section{Conclusion\label{sec5}}

This paper tried to model the idea that the decision-maker at time $t$ cannot
commit her successors by imagining that she can commit her immediate
successors, those in the interval $\left[  t,\ t+\varepsilon\right]  $, and
letting $\varepsilon\rightarrow0$. We then gave a rigourous definition of
(subgame perfect) equilibrium strategies, characterize them through the
equations (IE) and (DE). We give a local existence result in the\ analytic case.

One would, of course, like to have a global existence result, on $R^{d}%
\times\left[  0,\ T\right]  $, and to have weaker regularity assumptions
(sufficiently differentiable instead of analytic). Unfortunately, proving such
a theorem presents us with some serious mathematical challenges, and much more
work is required before we understand the situation. In fact, it seems to be
very similar to the situation wich prevailed on the (HJB)\ equation itself
before the discovery of viscosity solutions by Mike Crandall and Pierre-Louis
Lions when issues of existence, uniqueness and regularity where intertwined in
a very unsatisfactory manner. We feel that a similar program has to be
undertaken for equation (DE).

Another question is:\ why the continuous time ?\ Would it not be easier to
work with discrete time, and actually get the continous case by an appropriate
limiting process from the discrete case ?\ The answer, as we pointed out in
the introduction, is that we have no existence result for subgame perfect
equilibrium in the discrete case, so it is by no means clear that it is easier
than the continuous case. When Aumann started the study of economies with a
continuum of consumers, theorems were first proved directly, and the
connection with economies with a large number of consumers came much later.
For instance, the fact that if one constructs an economy with $Nn$ agents by
replicating $N$ times an economy with $n$ agents, the core of the large
economy converges to the equilibria of the limiting economy (which has a
continuum of agents) was first proved by Herbert Scarf, and was hailed as a
major achievement. Here again, such limiting results may hold for equilibrium
strategies, but it is another research program.

Finally, the obvious economic question is whether equilibrium strategies are
observationally different from optimal strategies. The difficulty here is
that, although the equilibrium strategy is defined for all $\left(
t,k\right)  $, we only observe one trajectory of the dynamics, the one that
starts at $k_{0}$ at time $t=0$. Devising testable consequences for our model
will be a third research program.

\appendix

\section{Proof of Theorem \ref{Th: condition necessaire}\textit{:}}

\subsection{Preliminaries}

Before proceeding with the proof of the theorem, let us mention some facts
about the flow $\mathcal{K}$ defined by (\ref{eq: HJBid}), (\ref{eq: HJBidic}).

Note first that the solution of (\ref{eq: HJBid}) which takes the value $k$ at
time $t$ coincides with the solution of the same equation which takes the
value $\mathcal{K}\left(  s_{1},t,k\right)  $ at time $s_{1}.$ In mathematical
terms, this property may be stated as
\begin{equation}
\mathcal{K}\left(  s_{2},t,k\right)  =\mathcal{K}\left(  s_{2},s_{1}%
,\mathcal{K}\left(  s_{1},t,k\right)  \right)  ,~~~0<t<s_{1}<s_{2}%
<T.\label{Rdf}%
\end{equation}

Next, consider the linearized equation around a prescribed solution
$t\rightarrow k_{0}\left(  t\right)  $ of the nonlinear system (\ref{25}), namely:%

\begin{equation}
\frac{dk}{ds}=\left(  \frac{\partial f}{\partial k}\left(  s,k_{0}\left(
s\right)  \right)  -\frac{\partial\sigma}{\partial k}\left(  s,k_{0}\left(
s\right)  \right)  \right)  k\left(  s\right) \label{a26}%
\end{equation}

This is a linear equation, so the flow is linear. The value at time $s$ of the
solution which takes the value $k$ at time $t\,\ $is $\mathcal{R}(s,t)k$,
where the matrix $\mathcal{R}(\cdot,t):[t,T]\longrightarrow R^{d\times d}$
satisfies:
\begin{align}
\frac{d\mathcal{R}}{ds}(s,t)  &  =\left(  \frac{\partial f}{\partial k}\left(
s,k_{0}\left(  s\right)  \right)  -\frac{\partial\sigma}{\partial k}\left(
s,k_{0}\left(  s\right)  \right)  \right)  \mathcal{R}\left(  s,t\right)
\label{a31}\\
\mathcal{R}\left(  t,t\right)   &  =\mathbf{I}\label{a32}%
\end{align}

From standard theory, it is well known that, if $f$ and $\sigma$ are $C^{k}$,
then $\mathcal{K}$ is $C^{k-1}$, and:
\begin{align*}
\frac{\partial\mathcal{K}\left(  s,t,k\right)  }{\partial k}  &
=\mathcal{R}(t,s)\\
\frac{\partial\mathcal{K}\left(  s,t,k\right)  }{\partial t}  &
=-\mathcal{R}(t,s)\left(  f\left(  t,k\right)  -\sigma\left(  t,k\right)
\right)
\end{align*}
where $\mathcal{R}(t,s)$ is computed by setting $k_{0}\left(  s\right)
=\mathcal{K}\left(  s,t,k\right)  $ in formulas (\ref{a31}), (\ref{a32}).

Let us now turn to the actual proof of Theorem \ref{Th: condition necessaire}:

\subsection{Necessary condition}

Given an equilibrium strategy $\sigma$, define a function $V$ by:
\begin{equation}
V(t,k)=\int_{t}^{T}h(s-t)u\left(  \sigma(s,\mathcal{K}\left(  s,t,k\right)
)\right)  ds+h(T-t)g\left(  \mathcal{K}\left(  T,t,k\right)  \right)
,\label{Vf}%
\end{equation}

Differentiating with respect to $k$, we find that:%
\begin{align*}
\frac{\partial V}{\partial k}\left(  t,k\right)   &  =\int_{t}^{T}%
h(s-t)\frac{\partial u}{\partial c}\left(  \sigma(s,\mathcal{K}\left(
s,t,k\right)  )\right)  \frac{\partial\sigma}{\partial k}\left(
s,\mathcal{K}\left(  s,t,k\right)  \right)  \frac{\partial\mathcal{K}%
}{\partial k}\left(  s,t,k\right)  ds\\
&  +h(T-t)\frac{\partial g}{\partial k}\left(  \mathcal{K}\left(
T,t,k\right)  \right)  \frac{\partial\mathcal{K}}{\partial k}\left(
T,t,k\right) \\
&  =\int_{t}^{T}h(s-t)\frac{\partial u}{\partial c}\left(  \sigma
(s,\mathcal{K}\left(  s,t,k\right)  )\right)  \frac{\partial\sigma}{\partial
k}\left(  s,\mathcal{K}\left(  s,t,k\right)  \right)  \mathcal{R}(s,t)ds\\
&  +h(T-t)\frac{\partial g}{\partial k}\left(  \mathcal{K}\left(
T,t,k\right)  \right)  \mathcal{R}(T,t)
\end{align*}

The (IE) equation will be derived by maximizing the individual payoff
$P_{1}\left(  t,k,\sigma,c\right)  $, as in formula (\ref{a14}). To this end,
let us first notice that the function $k_{1}$ defined by (\ref{RI}) and
(\ref{RCI}) can be written as
\[
k_{1}(s)=\mathcal{R}(s,t)\left(  \sigma(t,k)-c\right)  ,
\]
so that the individual payoff becomes
\begin{align*}
P_{1}\left(  t,k,\sigma,c\right)   &  =u\left(  c\right)  -u\left(
\sigma(t,k)\right) \\
&  +\int_{t}^{T}h\left(  s-t\right)  \frac{\partial u}{\partial c}\left(
\sigma\left(  s,k_{0}\left(  s\right)  \right)  \right)  \frac{\partial\sigma
}{\partial k}\left(  s,k_{0}\left(  s\right)  \right)  \mathcal{R}\left(
s,t\right)  \left(  \sigma(t,k)-c\right)  ds\\
&  +h\left(  T-t\right)  \frac{\partial g}{\partial k}\left(  k_{0}\left(
T\right)  \right)  \mathcal{R}(T,t)\left(  \sigma(t,k)-c\right)
\end{align*}

Since $u$ is concave and differentiable, the necessary and sufficient
condition to maximize $P_{1}(t,k,\sigma,c)$ with respect to $c$ is
\begin{align*}
\frac{\partial u}{\partial c}(c)  &  =\int_{t}^{T}h(s-t)\frac{\partial
u}{\partial c}\left(  \sigma(s,\mathcal{K}\left(  s,t,k\right)  )\right)
\frac{\partial\sigma}{\partial k}\left(  s,\mathcal{K}\left(  s,t,k\right)
\right)  \mathcal{R}(s,t)ds\\
&  +h(T-t)\frac{\partial g}{\partial k}\left(  \mathcal{K}\left(
T,t,k\right)  \right)  \mathcal{R}(T,t) = \frac{\partial V}{\partial k}(t,k).
\end{align*}

Therefore, the equilibrium strategy must satisfy
\[
\frac{\partial u}{\partial c}\left(  \sigma(t,k)\right)  =\frac{\partial
V}{\partial k}(t,k)
\]
and, substituting back into equation (\ref{Vf}), gives the (IE) equation.

\subsection{Sufficient condition}

Assume now that there exists a function $V$ satisfying (IE) and (BC), and
consider the strategy $\sigma=i\circ\frac{\partial V}{\partial k}$. Given any
consumption choice $c\in R^{d}$, the payoff to the decision-maker at time $t$
is:
\begin{align*}
P_{1}\left(  t,k,\sigma,c\right)   &  =u\left(  c\right)  -u\left(
\sigma(t,k)\right) \\
&  +\left[  \int_{t}^{T}h\left(  s-t\right)  \frac{\partial u}{\partial
c}\left(  \sigma\left(  s,k_{0}\left(  s\right)  \right)  \right)
\frac{\partial\sigma}{\partial k}\left(  s,k_{0}\left(  s\right)  \right)
\mathcal{R}(s,t)ds\right. \\
&  \left.  +h\left(  T-t\right)  \frac{\partial g}{\partial k}\left(
k_{0}\left(  T\right)  \right)  \mathcal{R}(T,t)\right]  \left(  \sigma\left(
t,k\right)  -c\right) \\
&  =u\left(  c\right)  -u\left(  \sigma(t,k)\right)  +\frac{\partial
V}{\partial k}(t,k)\left(  \sigma\left(  t,k\right)  -c\right) \\
&  =u\left(  c\right)  -u\left(  \sigma(t,k)\right)  -\frac{\partial
u}{\partial c}\left(  \sigma\left(  t,k\right)  \right)  \left(
c-\sigma\left(  t,k\right)  \right) \\
&  \leq0,
\end{align*}
where the first equality follows from the definition of $\mathcal{R}$, the
second equality is obtained by differentiating $V$ with respect to $k$, the
third equality follows from the definition of $\sigma$, and the last
inequality is due to the concavity of $u$. Observing that $P_{1}%
(t,k,\sigma,\sigma\left(  t,k\right)  )=0$, we see that the inequality
$P_{1}(t,k,\sigma,c)\leq0$ proves that $c=\sigma\left(  t,k\right)  $ achieves
the maximum so that $\sigma$ is an equilibrium strategy. Q.E.D.

\section{Proof of Proposition \ref{Proposition: IE / DE}}

Let a function $V:[0,T]\times R^{d}\rightarrow R^{d}$ be given. Consider the
function:%
\[
\varphi\left(  t,k\right)  =V\left(  t,k\right)  -\int_{t}^{T}h(s-t)u\left(
\sigma(s,\mathcal{K}\left(  s,t,k\right)  )\right)  ds-h(T-t)g\left(
\mathcal{K}\left(  T,t,k\right)  \right)
\]
where $\sigma=i\circ\frac{\partial V}{\partial k}$.

Consider the value of $\varphi$ along the trajectory of (\ref{dynamics})
originating from $k$ at time $t$. It is given by:%
\begin{align*}
\psi\left(  s,t,k\right)   &  =\varphi\left(  s,\mathcal{K}\left(
s,t,k\right)  \right) \\
&  =V\left(  s,\mathcal{K}\left(  s,t,k\right)  \right) \\
&  -\int_{s}^{T}h(x-s)u\left(  \sigma(x,\mathcal{K}\left(  x,s,\mathcal{K}%
\left(  s,t,k\right)  \right)  )\right)  dx\\
&  -h(T-s)g\left(  \mathcal{K}\left(  T,s,\mathcal{K}\left(  s,t,k\right)
\right)  \right) \\
&  =V\left(  s,\mathcal{K}\left(  s,t,k\right)  \right) \\
&  -\int_{s}^{T}h(x-s)u\left(  \sigma(x,\mathcal{K}\left(  x,t,k\right)
)\right)  dx\\
&  -h(T-s)g\left(  \mathcal{K}\left(  T,t,k\right)  \right)
\end{align*}
where we have used formula (\ref{Rdf}).

We compute the derivative of this function with respect to $s$:%
\begin{align*}
\frac{\partial\psi}{\partial s}\left(  s,k,t\right)   &  =\frac{\partial
V}{\partial t}\left(  s,\mathcal{K}\left(  s,t,k\right)  \right)
+\frac{\partial V}{\partial k}\left(  s,\mathcal{K}\left(  s,t,k\right)
\right)  \left(  f\left(  s,\mathcal{K}\left(  s,t,k\right)  \right)
-\sigma\left(  s,\mathcal{K}\left(  s,t,k\right)  \right)  \right) \\
&  +\int_{s}^{T}h^{\prime}(x-s)u\left(  \sigma(x,\mathcal{K}\left(
x,t,k\right)  )\right)  ds+u\left(  \sigma(s,\mathcal{K}\left(  s,t,k\right)
\right) \\
&  +h^{\prime}(T-s)g\left(  \mathcal{K}\left(  T,t,k\right)  \right)
\end{align*}

Since $\sigma=i\circ\frac{\partial V}{\partial k}$, we have:%
\[
-\frac{\partial V}{\partial k}\left(  s,\mathcal{K}\left(  s,t,k\right)
\right)  \sigma\left(  s,\mathcal{K}\left(  s,t,k\right)  \right)  +u\left(
\sigma(s,\mathcal{K}\left(  s,t,k\right)  \right)  =\tilde{u}\left(
\frac{\partial V}{\partial k}\left(  s,\mathcal{K}\left(  s,t,k\right)
\right)  \right)
\]

Substituting in the above, we get:%
\begin{align*}
\frac{\partial\psi}{\partial s}\left(  s,k,t\right)   &  =\frac{\partial
V}{\partial t}\left(  s,\mathcal{K}\left(  s,t,k\right)  \right)
+\frac{\partial V}{\partial k}\left(  s,\mathcal{K}\left(  s,t,k\right)
\right)  f\left(  s,\mathcal{K}\left(  s,t,k\right)  \right) \\
&  +\int_{s}^{T}h^{\prime}(x-s)u\left(  \sigma(x,\mathcal{K}\left(
x,t,k\right)  )\right)  ds+\tilde{u}\left(  \frac{\partial V}{\partial
k}\left(  s,\mathcal{K}\left(  s,t,k\right)  \right)  \right) \\
&  +h^{\prime}(T-s)g\left(  \mathcal{K}\left(  T,t,k\right)  \right)
\end{align*}

If (DE) holds, then $\psi\left(  s,k,t\right)  =$ $\psi\left(  T,k,t\right)
$, and if (BC) holds, then $\psi$ is identically zero, so that (IE) holds.
Conversely, if (IE) holds, then (BC) and (DE) obviously hold. Q.E.D.

\section{Proof of Proposition \label{G}\ref{phe}}

Define the function
\[
\mathcal{M}(s,t,k)=\mathcal{K}(s,t,k)e^{-\int_{t}^{s}(r(u)-\lambda(u))du}%
\]
with the understanding that the flow $\mathcal{K}$ is associated to the
strategy (\ref{csep}). Then we have%

\[
\frac{d\mathcal{M}\left(  s,t,k\right)  }{ds} =w(s) e^{- \int_{t}^{s} (r(u) -
\lambda(u) )du } - \lambda(s) e^{- \int_{t}^{s} (r(u) - \lambda(u) )du }
\int_{s}^{\infty} e^{-\int_{s}^{u} r(v) dv} w(u) du
\]
and after integrating this equation on $[t,s]$ we get
\begin{align}
\mathcal{M}\left(  s,t,k\right)  = k  &  + \int_{t}^{s} e^{-\int_{t}^{u} (
r(x) - \lambda(x) ) dx} w(u) du\nonumber\\
&  -\int_{t}^{s} \lambda(u) e^{- \int_{t}^{u} (r(x) - \lambda(x) )dx } \left[
\int_{u}^{\infty} e^{-\int_{u}^{v} r(x) dx} w(v) dv \right]
du.\label{eq: Appendix C1}%
\end{align}
The last term of this equality may be transformed into
\begin{align*}
&  \int_{t}^{s} \lambda(u) e^{- \int_{t}^{u} (r(x) - \lambda(x) )dx } \left[
\int_{u}^{\infty} e^{-\int_{u}^{v} r(x) dx} w(v) dv \right]  du\\
&  = \int_{t}^{s} \lambda(u) e^{\int_{t}^{u} \lambda(x) dx } \left[  \int
_{t}^{\infty} e^{-\int_{t}^{v} r(x) dx} w(v) dv - \int_{t}^{u} e^{-\int
_{t}^{v} r(x) dx} w(v) dv\right]  du\\
&  = \left[  \int_{t}^{\infty} e^{-\int_{t}^{u} r(x) dx} w(u) du \right]
\left[  \int_{t}^{s} \lambda(u) e^{\int_{t}^{u} \lambda(x) dx } du \right]  -
\int_{t}^{s} \lambda(u) e^{\int_{t}^{u} \lambda(x) dx } \left[  \int_{t}%
^{u}e^{-\int_{t}^{v} r(x) dx} w(v) dv \right]  du\\
&  = \left[  \int_{t}^{\infty} e^{-\int_{t}^{u} r(x) dx} w(u) du \right]
\left[  e^{\int_{t}^{s} \lambda(x) dx } -1\right]  - \int_{t}^{s} e^{-\int
_{t}^{v} r(x) dx} w(v) \left[  \int_{v}^{s} \lambda(u) e^{\int_{t}^{u}
\lambda(x) dx } du \right]  dv\\
&  = \left[  \int_{t}^{\infty} e^{-\int_{t}^{u} r(x) dx} w(u) du \right]
\left[  e^{\int_{t}^{s} \lambda(x) dx } -1\right]  - \int_{t}^{s} e^{-\int
_{t}^{v} r(x) dx} w(v) \left[  e^{\int_{t}^{s} \lambda(x) dx } - e^{\int
_{t}^{v} \lambda(x) dx }\right]  dv\\
&  = \left[  \int_{t}^{\infty} e^{-\int_{t}^{u} r(x) dx} w(u) du \right]
\left[  e^{\int_{t}^{s} \lambda(x) dx } -1\right]  - e^{\int_{t}^{s}
\lambda(x) dx } \int_{t}^{s} e^{-\int_{t}^{u} r(x) dx} w(u) du\\
&  ~~~~~~~~~~~~~~~+ \int_{t}^{s} e^{-\int_{t}^{u} \left(  r(x) - \lambda(x)
\right)  dx} w(u) du
\end{align*}
where the third equality follows from Fubini Theorem. Substituting this
formulation in equation $(\ref{eq: Appendix C1})$ yields
\[
\mathcal{M}\left(  s,t,k\right)  = k + \int_{t}^{\infty} e^{-\int_{t}^{u} r(x)
dx} w(u) du - e^{\int_{t}^{s} \lambda(x) dx } \left[  \int_{s}^{\infty}
e^{-\int_{t}^{u} r(x) dx} w(u) du \right]
\]
and therefore,
\[
\mathcal{K}\left(  s,t,k\right)  = e^{\int_{t}^{s} \left(  r(x) - \lambda(x)
\right)  dx} \left[  k + \int_{t}^{\infty} e^{-\int_{t}^{u} r(x) dx} w(u)
du\right]  - \int_{s}^{\infty} e^{-\int_{s}^{u} r(x) dx} w(u) du.
\]

Denoting by $V(t,k)$ the utility associated to the strategy $\sigma$ defined
by (\ref{csep}), we see that
\begin{align*}
V(t,k)  &  =\int_{t}^{\infty}\frac{\left[  \sigma(s,\mathcal{K}\left(
s,t,k\right)  )\right]  ^{1-\gamma}}{1-\gamma}h(s-t)ds\\
&  =\int_{t}^{\infty}\frac{\left[  \lambda(s)e^{\int_{t}^{s}(r(x)-\lambda
(x))dx}\right]  ^{1-\gamma}}{1-\gamma}\left[  k+\int_{t}^{\infty}e^{-\int
_{t}^{u}r(x)dx}w(u)du\right]  ^{1-\gamma}h(s-t)ds\\
&  =\frac{\left[  k+\int_{t}^{\infty}e^{-\int_{t}^{u}r(x)dx}w(u)du\right]
^{1-\gamma}}{1-\gamma}\int_{t}^{\infty}\left[  \lambda(s)e^{\int_{t}%
^{s}(r(x)-\lambda(x))dx}\right]  ^{1-\gamma}h(s-t)ds
\end{align*}
and that,
\[
\frac{\partial V}{\partial k}(t,k)=\left[  k+\int_{t}^{\infty}e^{-\int_{t}%
^{u}r(x)dx}w(u)du\right]  ^{-\gamma}\int_{t}^{\infty}\left[  \lambda
(s)e^{\int_{t}^{s}(r(x)-\lambda(x))dx}\right]  ^{1-\gamma}h(s-t)ds.
\]
Now, since the recursion $(\ref{Recursion})$ is satisfied, we see that
\[
\frac{\partial V}{\partial k}(t,k)=\left[  k+\int_{t}^{\infty}e^{-\int_{t}%
^{u}r(x)dx}w(u)du\right]  ^{-\gamma}\left[  \lambda(t)\right]  ^{-\gamma
}=\left(  \sigma(t,k)\right)  ^{-\gamma}.
\]
Therefore, the integrated equation $(\ref{IEcs})$ is satisfied which in turn
establishes that the strategy $\sigma$ defined by (\ref{csep}) is an
equilibrium strategy for the non-stationary problem $(\ref{Cs1}),(\ref{Cs3})$ .

\section{Proof of Theorem \ref{Theorem: indeterminacy}\label{H}}

\begin{proof}
Write equation (DE) for the function $V\left(  k\right)  $:%
\[
\tilde{u}\left(  V^{\prime}\right)  +fV^{\prime}+\int_{0}^{\infty}h^{\prime
}\left(  t\right)  u\left(  c\left(  t\right)  \right)  dt=0
\]
with $c\left(  t\right)  =\sigma\left(  \mathcal{K}\left(  t,k\right)
\right)  =i\left(  V^{\prime}\left(  \mathcal{K}\left(  t,k\right)  \right)
\right)  $.

Differentiate it at the equilibrium point $\bar{k}$. We get:%
\begin{equation}
\left(  \tilde{u}^{\prime}\left(  V^{\prime}\left(  \bar{k}\right)  \right)
+f\left(  \bar{k}\right)  \right)  V^{\prime\prime}\left(  \bar{k}\right)
+f^{\prime}\left(  \bar{k}\right)  V^{\prime}\left(  \bar{k}\right)  +\int
_{0}^{\infty}h^{\prime}\left(  t\right)  u^{\prime}\left(  \bar{c})\right)
i^{\prime}\left(  V^{\prime}\left(  \bar{k}\right)  \right)  V^{\prime\prime
}\left(  \bar{k}\right)  \frac{\partial\mathcal{K}\left(  t,k\right)
}{\partial k}dt\label{63}%
\end{equation}

We have $\tilde{u}^{\prime}\left(  V^{\prime}\left(  \bar{k}\right)  \right)
=-i\left(  V^{\prime}\left(  \bar{k}\right)  \right)  =-\bar{c}$, so that the
first term vanishes.

We are left with the two others. Note first that $i^{\prime}\left(  c\right)
=-\tilde{u}^{\prime\prime}\left(  c\right)  =-1/u^{\prime\prime}\left(
c\right)  $, so that the last integral can be rewritten as follows:%
\[
-\frac{u^{\prime}\left(  \bar{c}\right)  }{u^{\prime\prime}\left(  \bar
{c}\right)  }V^{\prime\prime}\left(  \bar{k}\right)  \int_{0}^{\infty
}h^{\prime}\left(  t\right)  \frac{\partial\mathcal{K}\left(  t,k\right)
}{\partial k}dt
\]

The function $y\left(  t\right)  =\partial\mathcal{K}\left(  t,k\right)
/\partial k$ is the solution of the linearized system at $\bar{k}$:%
\begin{equation}
\frac{dy}{dt}=\left(  f^{\prime}\left(  \bar{k}\right)  -i^{\prime}\left(
V^{\prime}\left(  \bar{k}\right)  \right)  V^{\prime\prime}\left(  \bar
{k}\right)  \right)  y=\left(  f^{\prime}\left(  \bar{k}\right)
-\alpha\right)  y\label{ekj}%
\end{equation}

Since $\bar{k}$ is an attractor, the exponent $\left(  f^{\prime}\left(
\bar{k}\right)  -\alpha\right)  $ must be non-positive, so $f^{\prime}\left(
\bar{k}\right)  -\alpha\leq0$, which is condition (\ref{59}).

From now on, we assume $f^{\prime}\left(  \bar{k}\right)  -\alpha<0\,$, so
that the linearized equation (\ref{ekj}) converges, and we have:
\[
\frac{\partial\mathcal{K}\left(  t,k\right)  }{\partial k}=\exp\left[  \left(
f^{\prime}\left(  \bar{k}\right)  -\alpha\right)  t\right]
\]

This gives us the last term in (\ref{63}). We now compute the middle term by
differentiating the formula for $V:$%
\[
V\left(  k\right)  =\int_{0}^{\infty}h\left(  t\right)  u\left(  c\left(
t\right)  \right)  dt
\]
yielding, by the same computation:%
\begin{align*}
V^{\prime}\left(  \bar{k}\right)   &  =-\frac{u^{\prime}\left(  \bar
{c}\right)  }{u^{\prime\prime}\left(  \bar{c}\right)  }V^{\prime\prime}\left(
\bar{k}\right)  \int_{0}^{\infty}h\left(  t\right)  \frac{\partial
\mathcal{K}\left(  t,k\right)  }{\partial k}dt\\
&  =-\frac{u^{\prime}\left(  \bar{c}\right)  }{u^{\prime\prime}\left(  \bar
{c}\right)  }V^{\prime\prime}\left(  \bar{k}\right)  \int_{0}^{\infty}h\left(
t\right)  \exp\left[  \left(  f^{\prime}\left(  \bar{k}\right)  -\alpha
\right)  t\right]  dt
\end{align*}

Substituting in equation (\ref{63}), we get:%
\[
\frac{u^{\prime}\left(  \bar{c}\right)  }{u^{\prime\prime}\left(  \bar
{c}\right)  }V^{\prime\prime}\left(  \bar{k}\right)  \left[  f^{\prime}\left(
\bar{k}\right)  \int_{0}^{\infty}h\left(  t\right)  \exp\left[  \left(
f^{\prime}\left(  \bar{k}\right)  -\alpha\right)  t\right]  dt+\int
_{0}^{\infty}h^{\prime}\left(  t\right)  \exp\left[  \left(  f^{\prime}\left(
\bar{k}\right)  -\alpha\right)  t\right]  dt\right]  =0
\]
and hence:%
\begin{align*}
f^{\prime}\left(  \bar{k}\right)   &  =-\frac{\int_{0}^{\infty}h^{\prime
}\left(  t\right)  \exp\left[  \left(  f^{\prime}\left(  \bar{k}\right)
-\alpha\right)  t\right]  dt}{\int_{0}^{\infty}h\left(  t\right)  \exp\left[
\left(  f^{\prime}\left(  \bar{k}\right)  -\alpha\right)  t\right]  dt}\\
&  =\frac{\left(  f^{\prime}\left(  \bar{k}\right)  -\alpha\right)  \int
_{0}^{\infty}h\left(  t\right)  \exp\left[  \left(  f^{\prime}\left(  \bar
{k}\right)  -\alpha\right)  t\right]  dt+1}{\int_{0}^{\infty}h\left(
t\right)  \exp\left[  \left(  f^{\prime}\left(  \bar{k}\right)  -\alpha
\right)  t\right]  dt}%
\end{align*}
where we have integrated by parts. This in turn gives formula (\ref{51}) and
concludes the proof.
\end{proof}

\section{Proof of Corollary \ref{corrolary: exponential discount}}

Substitute into equation (\ref{51}). We get for $\alpha$ the equation.%
\[
1=\alpha\int_{0}^{\infty}\exp\left[  -\rho+f^{\prime}\left(  \bar{k}\right)
-\alpha\right]  dt=\frac{\alpha}{\rho-f^{\prime}\left(  \bar{k}\right)
+\alpha}%
\]

If $f^{\prime}\left(  \bar{k}\right)  \neq\rho$, there is no solution to this
equation, so $f^{\prime}\left(  \bar{k}\right)  $ must be equal to $\rho$. In
order to establish (\ref{53bis}), we differentiate twice the (DE) equation
(which is here the HJB equation) and evaluate it at $\bar{k}$ and get
\[
-V^{\prime}(\bar{k})f^{\prime\prime}(\bar{k})=u^{\prime\prime}(\bar{c}%
)\alpha\left(  f^{\prime}(\bar{k})-\alpha\right)  .
\]
This is a quadratic equation in $\alpha$ and, assuming $f$ is concave, it
admits two roots,
\[
\alpha=f^{\prime}(\bar{k})\left[  \frac{1\pm\sqrt{1+4\frac{u^{\prime}(\bar
{k})f^{\prime\prime}(\bar{k})}{\rho^{2}u^{\prime\prime}(\bar{k})}}}{2}\right]
\]
The second root is not valid because $\bar{k}$ is an attractor and hence
$0<f^{\prime}(\bar{k})\leq\alpha$. This leave us with the only possible root
$\alpha$ given by (\ref{53bis}). Q.E.D.

\section{Proof of Proposition \ref{pol}}

Substituting the specification of the discount function into equation
(\ref{51}) gives
\begin{equation}
1=\alpha\int_{0}^{T}\exp\left[  \left(  -\rho+f^{\prime}\left(  \bar{k}
\right)  -\alpha\right)  t \right]  dt\label{52}%
\end{equation}

Since $\alpha\geq f^{\prime}\left(  \bar{k}\right)  $, the term $\rho
-f^{\prime}\left(  \bar{k}\right)  +\alpha$ is different from $0$ and
therefore the equation (\ref{52}) becomes
\begin{equation}
\label{eq: indeterminacy 1}\left(  f^{\prime}\left(  \bar{k}\right)
-\rho\right)  T \exp\left[  -T\left(  f^{\prime}\left(  \bar{k}\right)  -
\rho\right)  \right]  =\alpha T\exp\left(  -\alpha T\right)
\end{equation}

The left inequality of (\ref{55}) is simply obtained by noticing that, due to
equation (\ref{52}), $\alpha$ must be positive and the equation
(\ref{eq: indeterminacy 1}) implies then that $f^{\prime}\left(  \bar
{k}\right)  -\rho> 0$.

Now notice that equation (\ref{eq: indeterminacy 1}) is of the type
$xe^{x}=ye^{y}$. It has the obvious solution $x=y$, plus another one. The
first solution gives $f^{\prime}\left(  \bar{k}\right)  -\rho=\alpha$,
contradicting the fact that $\rho-f^{\prime}\left(  \bar{k}\right)
+\alpha\neq0$, so it must be rejected. The second solution defines $y$~as a
function of $x$, say $y=\varphi\left(  x\right)  $ for $x>0$, which is easily
seen to be decreasing and to obey the properties given in
(\ref{eq: varphi properties}). So $\alpha T=\varphi\left(  T\left(  f^{\prime
}\left(  \bar{k}\right)  -\rho\right)  \right)  $, and formula (\ref{62}) follows.

Now, taking into account the condition $\alpha\geq f^{\prime}\left(  \bar
{k}\right)  $ given by (\ref{59}), equation (\ref{62}) implies
\[
\varphi\left(  f^{\prime}\left(  \bar{k}\right)  T-\rho T\right)  \geq
f^{\prime}\left(  \bar{k}\right)  T > f^{\prime}\left(  \bar{k}\right)  T-\rho
T.
\]

Since the function $\varphi$ is decreasing and $\varphi(1) =1$, the above
inequality implies
\[
f^{\prime}\left(  \bar{k}\right)  T- \rho T <1
\]
which is precisely the right inequality of of (\ref{55}).

\end{document}